\documentclass[11pt]{article}
\usepackage[top=1in, bottom=1in, left=1in, right=1in]{geometry}
\usepackage{tikz}
\usetikzlibrary{shapes,arrows,matrix,patterns,positioning}
\usepackage{color,diagbox}
\usepackage{graphicx}
\usepackage{algorithmic}
\usepackage{amssymb}
\usepackage{epstopdf}
\usepackage[ruled]{algorithm2e}
\usepackage{float}
\usepackage{amsmath,amsthm,bm,color,epsfig,enumerate,caption}
\usepackage{hyperref}
\usepackage{booktabs}
\usepackage{multirow}
\usepackage{array}
\usepackage{tabu}
\usepackage{breqn}
\usepackage{mathtools}
\usepackage{multicol}
\usepackage{multirow}
\DeclarePairedDelimiter{\ceil}{\lceil}{\rceil}

\newtheorem{theorem}{Theorem}[section]

\newtheorem{algo}[theorem]{Algorithm}

\renewcommand{\appendix}[1]{
\section*{Appendix: #1}
}

\newcommand{\bbC}{\mathbb{C}}
\newcommand{\bbR}{\mathbb{R}}
\newcommand{\ii}{\mathfrak{i}}

\makeatletter
\newcommand*{\extendadd}{
  \mathbin{
    \mathpalette\extend@add{}
  }
}
\newcommand*{\extend@add}[2]{
  \ooalign{
    $\m@th#1\leftrightarrow$%
    \vphantom{$\m@th#1\updownarrow$}
    \cr
    \hfil$\m@th#1\updownarrow$\hfil
  }
}
\makeatother

\begin{document}

\title{Fast Algorithms for the Multi-dimensional\\ Jacobi Polynomial Transform}
\author{James Bremer \\ Department of Mathematics\\ University of California, Davis, CA, USA\\ \href{mailto:bremer@math.ucdavis.edu}{bremer@math.ucdavis.edu}
\and Qiyuan Pang \\ Department of Mathematics\\ Purdue University, West Lafayette, IN, USA\\  \href{mailto:qpang@purdue.edu}{qpang@purdue.edu} 
     \and Haizhao Yang \\ Department of Mathematics\\ Purdue University, West Lafayette, IN, USA\footnote{Current address.}\\ National University of Singapore, Singapore\footnote{On leave from the National University of Singapore, where part of the work was done.}\\ \href{mailto:haizhao@nus.edu.sg}{haizhao@nus.edu.sg} }
     
\date{\today}
\maketitle

\begin{abstract}

We use the well-known observation that the solutions of Jacobi's differential equation can be represented via the non-oscillatory phase and amplitude functions to develop 
 a fast algorithm for computing multi-dimensional Jacobi polynomial transforms.
More explicitly,  it  follows from this observation that the matrix corresponding to the 
discrete Jacobi transform is the Hadamard product of a 
numerically low-rank matrix and a multi-dimensional discrete Fourier transform (DFT) matrix.
The application of the Hadamard product can be carried out via $r^d$ fast Fourier transforms (FFTs), where $r=O(\dfrac{\log n}{\log \log n})$ and $d$ is the dimension, resulting in a nearly optimal algorithm to compute the multidimensional Jacobi polynomial transform. 

\end{abstract}


\section{Introduction}
\label{sec:intro}

{{Jacobi polynomials contain a wide class of orthogonal polynomials, e.g., the Gegenbauer polynomials, and thus also the Legendre, Zernike, and Chebyshev polynomials, which have been studied and applied extensively in mathematical analysis and practical applications \cite{doi:10.1137/1.9781611970470,szeg1939orthogonal,timan1994theory}. Previously, the Legendre and Chebyshev polynomials, not only in the univariate cases but also in the multivariate cases, have been essential tools in solving partial differential equations with spectral methods \cite{boyd2001chebyshev,Canuto2006,funaro1992polynomial,doi:10.1137/1.9781611970425,doi:10.1142/3662}. Fast Chebyshev polynomial transforms in 1D to 3D have been proposed in \cite{chebfun1d,chebfun2d,chebfun3d} to increase the numerical efficiency of these transforms when the problem size is large. Recently, Jacobi polynomials in spectral approximations were revisited and shown to be an efficient tool for problems with degenerated or singular coefficients \cite{BERNARDI1997209}, for singular differential equations  \cite{BENYU1998180,GUO2000373,10.2307/43692835}, and for optimal error estimates for $p$-version of finite element methods \cite{Guo,doi:10.1137/S0036142901356551,GUO2013122}.  The Jacobi polynomials were also generalized to a wider class of functions to simplify the approximation analysis with more precise error estimates and well-conditioned algorithms in \cite{GUO20091011}. Therefore, it is important to develop efficient computational tools for multivariate Jacobi polynomial transforms.}}

The one-dimensional forward discrete Jacobi transform consists of evaluating an expansion of the form
\begin{equation}
\label{eq:jactrs}
f(x) = \sum\limits_{\nu=1}^{n}\alpha_{\nu}P_{\nu-1}^{(a,b)}(x),
\end{equation}   
where $P_{\nu}^{(a,b)}$ denotes the order-$\nu$ Jacobi polynomial of the first kind corresponding to parameters $a$ and $b$, 
at a collection of points 
{{$X = \{x_i\}_{i=1,\ldots,n}\subset (-1,1)$}}.  
The one-dimensional inverse discrete Jacobi transform is the process of 
 computing the coefficients 
$\{\alpha_\nu\}_{i=1,\ldots,n}$
 in an expansion of the form (\ref{eq:jactrs}) given its values at a collection of 
distinct points 
{{$X$}}.    For the sake of brevity, we will generally drop the adjective ``discrete'' and simply use the
terms  ``forward Jacobi transform'' and  ``inverse Jacobi transform'' when referring to these operations.
Often, the points {{$X$}}  are the nodes of the $n$-point Gauss-Jacobi quadrature rule

\begin{equation}
\label{eq:GJquadrule}
\int_{-1}^{1}p(x)(1+x)^{a}(1+x)^{b}dx \approx \sum\limits_{\nu=1}^{n}p(x_{\nu})\omega_{\nu}
\end{equation}
that is exact when $p$ is a  polynomial of degree less than or equal to $2n-1$.  However, it is useful 
to consider  more general sets of points as well.  In the case in which the points 
{{$X$}} are the nodes
of the Gauss-Jacobi quadrature rule, we say the corresponding transform and inverse transform are uniform; otherwise, we describe them as nonuniform. To form an orthonormal matrix to represent the transform for numerical purpose, we will rescale the transform in \eqref{eq:jactrs} with quadrature weights later.

Many methods for rapidly applying the  Jacobi transform and various special cases
of the Jacobi transform have been developed.  Examples include the algorithms of 
 \cite{Iserles2011,hale2014fast,alpert1991fast,Rokhlin1} for applying the Legendre transform, and 
those of \cite{Keiner2} for forming  and evaluating expansions  for Jacobi polynomials in general.
See also,  \cite{chebfun1d} and its references for extensive information on numerical algorithms for forming and manipulating Chebyshev expansions.
Almost all such algorithms can be placed into one of two categories. Algorithms in the first category, such as
{{\cite{askey1975orthogonal,andrews1999special,Maroni2008,hua1963harmonic,szeg1939orthogonal,91376860,61979980,206962923,1998MaCom..67.1577P}}}, 
make use of the structure  of the connection matrices which take the coefficients in the expansion of a function in terms of one set of 
Jacobi polynomials to the coefficients in the expansion of the same function for a different set of Jacobi polynomials.
They typically operate by  computing the Chebyshev coefficients of expansion and then applying a connection matrix or a series
of connection matrices to obtain the coefficients in the desired expansion.
The computation of the Chebyshev expansion can be carried out efficiently in $O(n\log n)$ operations via the nonuniform FFT (NUFFT) \cite{Greengard,Townsend}
and the application of each connection matrix can be performed in a number of operations which grows
linearly or quasi-linearly 
via the fast multipole method (FMM) \cite{alpert1991fast,greengard1987fast,keiner2009computing,Keiner2} 
or the fast eigendecomposition of semiseparable matrices \cite{gemignani2010matrix,keiner2008gegenbauer,Keiner2}.

The second class of algorithms, of which \cite{Oneil-Rokhlin} is a major example,
make use of butterfly methods and similar techniques for the application of oscillatory
matrices.   Most algorithms in this category, including that of in \cite{Oneil-Rokhlin}, require precomputation with $O(n^2)$   running times.
However, \cite{Jacobi} describes an algorithm based on this approach whose total running time is $O(n \log^2(n))$, when both parameters $a$ and $b$ are in the interval $(-1/2,1/2)$.
It makes use of  the observation that the solutions of  Jacobi's differential equation can be accurately represented via 
non-oscillatory phase and amplitude functions to apply the Jacobi transform rapidly.
The fact that certain differential equations admit non-oscillatory phase functions has long been known \cite{DLMF}; a general theorem
 was established in \cite{HEITMAN20151} 
and a numerical method for the computation of non-oscillatory phase functions was 
given in \cite{BremerKummer}.
In the case of Jacobi's differential equation with {{$t\in [0,\pi]$ and $\nu > 0$}},
there exists a smooth amplitude function $M^{(a,b)}(t,\nu)$ and a smooth phase function $\psi^{(a,b)}(t,\nu)$ such that
\begin{equation}
\label{eq:jacnonosci}
\tilde{P}_{\nu}^{(a,b)}(t) = M^{(a,b)}(t,\nu)\cos(\psi^{(a,b)}(t,\nu))
\end{equation}  
and
\begin{equation}
\label{eq:jacnonosci2}
\tilde{Q}_{\nu}^{(a,b)}(t) = M^{(a,b)}(t,\nu)\sin(\psi^{(a,b)}(t,\nu)),
\end{equation}  
where $\tilde{P}_{\nu}^{(a,b)}$ and $\tilde{Q}_{\nu}^{(a,b)}$ are referred to as the modified Jacobi functions of the first and second kind, respectively.
They are defined in terms of the Jacobi functions of the first and second kinds $P_\nu^{(a,b)}$ and  $Q_\nu^{(a,b)}$ (see \cite{DLMF} for definitions)
via the formulas
\begin{equation}
\label{eq:modifiedjac}
\tilde{P}_{\nu}^{(a,b)}(t) = C_{\nu}^{(a,b)}P_{\nu}^{(a,b)}(\cos(t))\sin\left(\frac{t}{2}\right)^{a+\frac{1}{2}}\cos\left(\frac{t}{2}\right)^{b+\frac{1}{2}}
\end{equation}
and
\begin{equation}
\tilde{Q}_\nu^{(a,b)}(t) = 
C_{\nu}^{(a,b)} \ Q_\nu^{(a,b)}\left(\cos(t)\right) 
\sin\left(\frac{t}{2}\right)^{a+\frac{1}{2}}
\cos\left(\frac{t}{2}\right)^{b+\frac{1}{2}},
\label{eq:modifiedjac2}
\end{equation}
where
\begin{equation}
\label{eq:C}
C_{\nu}^{(a,b)} = \sqrt{(2\nu+a+b+1)\frac{\Gamma(1+\nu)\Gamma(1+\nu+a+b)}{\Gamma(1+\nu+a)\Gamma(1+\nu+b)}}.
\end{equation}
The normalization constant is chosen to ensure the $L^{2}(0,\pi)$ norm of $\tilde{P}_{\nu}^{(a,b)}$ is $1$ when $\nu$ is an integer. Indeed, the set $\{\tilde{P}_{j}^{(a,b)}\}_{j=0}^{\infty}$ 
 is an orthonormal basis for $L^{2}(0,\pi)$. The change of variables $x = \cos(t)$ makes the singularities in phase and amplitude functions for Jacobi's differential equation more tractable. 
Obviously, there is no substantial difference in treating expansions of the type (\ref{eq:jactrs})
 and those of the form
\begin{equation}
\label{eq:jactrs2}
f(t) = \sum\limits_{\nu=1}^{n}\alpha_{\nu}\tilde{P}_{\nu-1}^{(a,b)}(t)=\sum\limits_{\nu=1}^{n}\alpha_{\nu}M^{(a,b)}(t,\nu)\cos (\psi^{(a,b)}(t,\nu)),
\end{equation} 
and we will consider the latter form here, more specially, the latter form with scaling by quadrature weights.  Moreover, we will denote by $\{t_k\}_{k=1,\ldots,n}$ and 
$\{w_k\}_{k=1\ldots, n}$ the nodes and weights of the trigonometric Gauss-Jacobi quadrature rule

\begin{equation}
\int_{0}^\pi f(\cos(t)) \cos^{2a+1}\left(\frac{t}{2}\right)   \sin^{2b+1}\left(\frac{t}{2}\right)\ dt 
\approx
\sum_{k=1}^n f(\cos(t_k)) 
\cos^{2a+1}\left(\frac{t_k}{2}\right)   \sin^{2b+1}\left(\frac{t_k}{2}\right)
w_k
\label{introduction:modrule}
\end{equation}
obtained by applying the change of variables $x=\cos(t)$ to (\ref{eq:GJquadrule}).
Again, we refer to transforms which use the set of points $\{t_k\}_{k=1,\ldots,n}$ as uniform
and those which use some other point sets as nonuniform.






In \cite{Jacobi}, it was observed that the second summation in
\eqref{eq:jactrs2} could be interpreted as the application of a 
 structured matrix that is the real part of a Hadamard product of a NUFFT matrix and a numerically low-rank matrix.  This leads to a method for
its computation which takes quasi-linear time 
\cite{Townsend,Candes-Demanet-Ying,yang2018unified}. Expansions in terms of the functions of the second
kind  $\{\tilde{Q}_{j}^{(a,b)}\}_{j=0}^{\infty}$ can be handled similarly. This non-oscillatory representation method will be described in more detail in
Section \ref{sec:NPF}.

In this paper, we generalize the one-dimensional Jacobi polynomial transform of \cite{Jacobi} to the multi-dimensional case.
The transformation matrix in the multi-dimensional case is still the real part of a Hadamard product of a multi-dimensional NUFFT matrix and a numerically low-rank matrix in the form of tensor products. However, the numerical rank of the low-rank matrix might increase quickly in the dimension following the algorithm in \cite{Jacobi} and hence the multi-dimensional NUFFT in \cite{Greengard,Townsend} might not be sufficiently efficient. To obtain a viable multi-dimensional Jacobi polynomial transform, we reformulate the Hadamard product into a new one, which is the Hadamard product of an over-sampled Fourier transform matrix and a numerically low-rank matrix, and propose a fast randomized SVD to achieve near optimality in the low-rank approximation. This leads to an efficient implementation of two-dimensional and three-dimensional Jacobi polynomial transforms.     
Besides, as will be demonstrated numerically, the new method proposed in this paper is {{faster and }}more robust than the method in \cite{Jacobi} and works for $a$ and $b$ in the regime $(-1,1)$, {{which is a larger interval than}} the regime $(-\frac{1}{2},\frac{1}{2})$ in \cite{Jacobi}. 

The remainder of the paper is organized as follows. In Section \ref{sec:EJT}, we will first briefly introduce the fast SVD via randomized sampling and the theory of non-oscillatory phase functions. In Section \ref{sec:jt}, a variant of the fast one-dimensional Jacobi polynomial transform in \cite{Jacobi} is given; it is followed
by a description of our fast multi-dimensional Jacobi polynomial transform. In Section \ref{sec:results}, numerical results which demonstrate the efficiency of our algorithm are described. Finally, we will conclude our discussion in Section \ref{sec:conclusion}.

\section{Preliminaries}

\label{sec:EJT}

In this section, we will revisit the linear scaling randomized SVD introduced in \cite{engquist2009} and the non-oscillatory phase and amplitude functions in \cite{Jacobi} to make the presentation self-contained.

\subsection{Approximate SVD via Randomized Sampling}
\label{sec:SVD}
For a numerically low-rank matrix $Z\in\mathbb{C}^{n\times n}$ with an $O(1)$ algorithm to evaluate an arbitrary entry, \cite{engquist2009} introduced an $O(n r^2)$ algorithm to construct a rank-$r$ approximate SVD, $Z \approx U_{0}\Sigma_{0}V_{0}^{*}$, from $O(r)$ randomly selected rows and columns of $Z$, where $^*$ denotes the conjugate transpose, $U_0$ and $V_0\in\mathbb{C}^{n\times r}$, $\Sigma_{0}\in \mathbb{C}^{r\times r}$ is a diagonal matrix with approximate singular values in the diagonal part in a descent order.

Here, we adopt the standard notation for a submatrix in MATLAB: given a row index set $I$ and a column index set $J$, $Z_{I,J} = Z(I,J)$ is the submatrix with entries from rows in $I$ and columns in $J$; we also use ``$:$" to denote the entire columns or rows of the matrix, i.e., $Z_{I,:} = Z(I,:)$ and $Z_{:,J} = Z(:,J)$. With these handy notations, we briefly introduce the randomized SVD as follows.

\begin{algo}{Approximate SVD via randomized sampling}
\begin{enumerate}
\label{algo}
\item 
Let $\Pi_{col}$ and $\Pi_{row}$ denote the important columns and rows of $Z$ that are used to form the column and row bases. Initially $\Pi_{col} = \emptyset$ and $\Pi_{row} = \emptyset$.
\item
Randomly sample $rq$ rows and denote their indices by $S_{row}$. Let $I = S_{row}\cup \Pi_{row}$. Here $q = O(1)$ is a multiplicative oversampling parameter. Perform a pivoted $QR$ decomposition of $Z_{I,:}$ to get 
\begin{eqnarray}
Z_{I,:}P = QR,
\end{eqnarray}
where $P$ is the resulting permutation matrix and $R = (r_{ij})$ is an $O(r) \times n$ upper triangular matrix. Define the important column index set $\Pi_{col}$ to be the first $r$ columns picked within the pivoted $QR$ decomposition.
\item
 Randomly sample $rq$ columns and denote their indices by $S_{col}$. Let $J = S_{col}\cup \Pi_{col}$. Perform a pivoted $LQ$ decomposition of $Z_{:,J}$ to get
\begin{eqnarray}
PZ_{:,J} = LQ,
\end{eqnarray}
where $P$ is the resulting permutation matrix and $L = (l_{ij})$ is an $m \times O(r)$ lower triangular matrix. Define the important row index set $\Pi_{row}$ to be the first $r$ rows picked within the pivoted $LQ$ decomposition.
\item
Repeat steps 2 and 3 a few times to ensure $\Pi_{col}$ and $\Pi_{row}$ sufficiently sample the important columns and rows of $Z$.
\item
Apply the pivoted $QR$ factorization to $Z_{:.\Pi_{col}}$ and let $Q_{col}$ be the matrix of the first $r$ columns of the $Q$ matrix. Similarly, apply the pivoted $QR$ factorization to $Z_{\Pi_{row},:}^{*}$ and let $Q_{row}$ be the
matrix of the first $r$ columns of the $Q$ matrix.
\item
We seek a middle matrix $M$ such that $Z \approx Q_{col}MQ_{row}^{*}$. To solve this problem efficiently, we approximately reduce it to a least-squares problem of a smaller size. Let $S_{col}$ and $S_{row}$ be the index sets of a few extra randomly sampled columns and rows. Let $J = \Pi_{col} \cup S_{col}$ and
$I = \Pi_{row} \cup S_{row}$. A simple least-squares solution to the problem
\begin{eqnarray}
\min\limits_{M}\parallel Z_{I,J}-(Q_{col})_{I,:}M(Q_{row}^{*})_{:,J}\parallel,
\end{eqnarray}
gives $M = (Q_{col})_{I,:}^{\dagger}Z_{I,J}(Q_{row}^{*})_{:,J}^{\dagger}$, where $(\cdot)^{\dagger}$ stands for the pseudo-inverse.
\item
Compute an SVD $M \approx U_{M}\Sigma_{M}V_{M}^{*}$. Then the low-rank approximation of $Z \approx U_{0}S_{0}V_{0}^{*}$ is given by
\begin{eqnarray}
U_{0} = Q_{col}U_{M}; \Sigma_{0} = \Sigma_{M}; V_{0}^{*} = V_{M}^{*} Q_{row}^{*}.
\end{eqnarray}
\end{enumerate}
\end{algo}

In our numerical implementation, iterating Steps $2$ and $3$ twice is empirically sufficient to achieve accurate low-rank approximations via Algorithm \ref{algo}. Similar arguments as in \cite{subCUR} for a randomized CUR factorization can be applied to quantify the error and success probability rigorously for Algorithm \ref{algo}. But at this point, we only focus on the application of Algorithm \ref{algo} to the fast Jacobi polynomial transform without theoretical analysis.

Note that our goal in the Jacobi polynomial transform is to construct a low-rank approximation of the low-rank matrix, i.e., $Z\approx U V^*$ with $U$ and $V\in\mathbb{C}^{n\times r}$, up to a fixed relative error $\varepsilon$, rather than a fixed rank. Algorithm \ref{algo} can also be embedded into an iterative process that gradually increases the rank parameter $r$ to achieve the desired accuracy. We denote the rank parameter as $r_\epsilon$ when it is large enough to obtain the accuracy $\epsilon$.

When $Z$ origins from the discretization of a smooth function $Z(x,y)$ at the grid points $\{x_i\}_{1\leq i\leq n}$ and $\{y_i\}_{1\leq i\leq n}$, i.e., $Z_{ij} = Z(x_i,y_j)$, another standard method for constructing a low-rank factorization of $Z\approx UV^*$ is Lagrange interpolation at Chebyshev grids in $x$ or $y$. For example, let $C_\mu:=\{\mu_i\}_{1\leq i\leq r}$ denote the set of $r$ Chebyshev grids in the domain of $x$, $U\in \mathbb{R}^{n\times r}$ be the matrix consisting of the $i$-th Lagrange polynomial in $x$ corresponding to $\mu_i$ as its $i$-th column, and $V\in\mathbb{C}^{n\times r}$ be the matrix such that $V_{ij}$ is equal to $\bar{Z}(\mu_j,y_i)$, where $\bar{\cdot}$ denotes the conjugate operator, then $Z\approx UV^*$ by the Lagrange interpolation. Usually, an oversampling parameter $q$ is used via setting $r=q r_\epsilon$. Then a rank-$r_\epsilon$ truncated SVD of $U\approx U_0 \Sigma_0 V_0^*$ gives a compressed rank-$r_\epsilon$ approximation of $Z\approx U_0 \left( V V_0 \Sigma_0\right)^*$, where $U_0\in\mathbb{C}^{n\times r_\epsilon}$ and $V V_0 \Sigma_0\in\mathbb{C}^{n\times r_\epsilon}$.

The fast nonuniform FFT in \cite{Townsend} and the one-dimensional Jacobi polynomial transform in \cite{Jacobi} adopted low-rank approximation via Lagrange interpolation to deal with the low-rank term in their Hadamard products of low-rank and (nonuniform) FFT matrices without an extra truncated SVD. In this paper, we propose to use the randomized SVD via random sampling to obtain nearly optimal rank in the low-rank approximation.

\subsection{Non-oscillatory phase and amplitude functions}
\label{sec:NPF}
Given a pair of parameters $a$ and $b$ in $\left(-\frac{1}{2},\frac{1}{2}\right)$, and a 
maximum degree of interest $N_{\mbox{\tiny max}} > 27$, we revisit the fast algorithms in \cite{Jacobi} for constructing non-oscillatory phase and amplitude functions $\psi^{(a,b)}(t,\nu)$ and $M^{(a,b)}(t,\nu)$ such that 
\begin{equation}
\tilde{P}_\nu^{(a,b)}(t) = M^{(a,b)}(t,\nu) \cos\left(\psi^{(a,b)}(t,\nu)\right)
\label{introduction:p}
\end{equation}
and
\begin{equation}
\tilde{Q}_\nu^{(a,b)}(t) = M^{(a,b)}(t,\nu) \sin\left(\psi^{(a,b)}(t,\nu)\right)
\label{introduction:q}
\end{equation}
for $t\in[\frac{1}{N_{\mbox{\tiny max}}},\pi-\frac{1}{N_{\mbox{\tiny max}}}]$ and $\nu\in(27,N_{\mbox{\tiny max}}]$. The polynomials with $t$ and $\nu$ out of these ranges can be evaluated by the well-known three-term recurrence relations or various asymptotic expansions; the lower bound of $27$ was chosen to obtain optimal numerical performance.

Next, we present a few facts regarding the phase and amplitude functions related to  Jacobi's differential equations. It is well known 
that the functions $\tilde{P}_\nu^{(a,b)}$ and $\tilde{Q}_\nu^{(a,b)}$ satisfy the second order differential equation
\begin{equation}
y''(t) + q^{(a,b)}_\nu(t) y(t) = 0,
\label{jacobieq:mod}
\end{equation}
where
\begin{dmath}
q_\nu^{(\alpha,\beta)}(t) = 
\left(\nu + \frac{\alpha+\beta+1}{2} \right)^2 
+ \frac{\frac{1}{4}-\alpha^2}{4 \sin\left(\frac{t}{2}\right)^2}
+ \frac{\frac{1}{4}-\beta^2}{4 \cos\left(\frac{t}{2}\right)^2}.
\label{jacobieq:coef}
\end{dmath}

We refer to Equation (\ref{jacobieq:mod}) as Jacobi differential
equation. Following the derivation in \cite{Jacobi}, we can show that the pair $\{\tilde{P}_\nu^{(a,b)},\tilde{Q}_\nu^{(a,b)}\}$ of real-valued,  linearly independent solutions of \eqref{jacobieq:mod} satisfies 
\begin{equation}
\tilde{P}_\nu^{(a,b)}(t) = \sqrt{W}\ \frac{\cos\left(\psi^{(a,b)}(t,\nu)\right)}{\sqrt{\left|\partial_t \psi^{(a,b)}(t,\nu)\right|}},
\label{phase:u2}
\end{equation}
and
\begin{equation}
\tilde{Q}_\nu^{(a,b)}(t) = \sqrt{W}\ \frac{\sin\left(\psi^{(a,b)}(t,\nu)\right)}{\sqrt{\left|\partial_t \psi^{(a,b)}(t,\nu)\right|}},
\label{phase:v2}
\end{equation}
where $W$ is the necessarily positive constant Wronskian  of the pair $\{\tilde{P}_\nu^{(a,b)}, \tilde{Q}_\nu^{(a,b)}\}$, and $\psi^{(a,b)}(t,\nu)$ is the non-oscillatory phase function of the pair $\{\tilde{P}_\nu^{(a,b)}, \tilde{Q}_\nu^{(a,b)}\}$ satisfying
\begin{equation}
\psi^{(a,b)}(t,\nu) = C + \int_{\sigma_1}^t 
 \frac{W}{(\tilde{P}_\nu^{(a,b)}(s))^2+(\tilde{Q}_\nu^{(a,b)}(s))^2} \ ds
\end{equation}
with $C$ an appropriately chosen constant related to $\sigma_1$. Note that $\partial_t\psi^{(a,b)}(t,\nu)>0$ since $W>0$. Hence, the non-oscillatory amplitude function of $\{\tilde{P}_\nu^{(a,b)}, \tilde{Q}_\nu^{(a,b)}\}$ can be defined as
\begin{equation}
M^{(a,b)}(t,\nu) = \sqrt{\frac{W}{\left|\partial_t \psi^{(a,b)}(t,\nu)\right|}}=\sqrt{\frac{W}{\partial_t \psi^{(a,b)}(t,\nu)}}.
\label{phase:mpsip}
\end{equation}
Through straightforward computation, 
it can be verified that
the square  $N^{(a,b)}(t,\nu) = \left(M^{(a,b)}(t,\nu)\right)^2$ of the amplitude function
satisfies the third order linear ordinary differential equation (ODE) 
\begin{equation}
\partial_{ttt} N^{(a,b)}(t,\nu) + 4 q_\nu^{(a,b)}(t)  \partial_t N^{(a,b)}(t,\nu)  + 2 (q_\nu^{(a,b)}(t))' N^{(a,b)}(t,\nu) = 0
\ \ \mbox{for all}\ \ \sigma_1 < t < \sigma_2.
\label{phase:Neq}
\end{equation}

Obviously,
\begin{equation}
\left(M^{(a,b)}(t,\nu)\right)^2 = 
\left(P_\nu^{(a,b)}(t)\right)^2+\left(Q_\nu^{(a,b)}(t)\right)^2.
\label{jacobiphase:M}
\end{equation}
Hence, {{to set up the initial condition for the ODE in \eqref{jacobiphase:M},}} we can specify the values of $\left(M^{(a,b)}(t,\nu)\right)^2$ and its first two derivatives in $t$ 
at a point on the interval $[\frac{1}{N_{\mbox{\tiny max}}},\pi-\frac{1}{N_{\mbox{\tiny max}}}]$ for any $\nu$ using various asymptotic expansion for 
 $P_\gamma^{(a,b)}$ and $Q_\gamma^{(a,b)}$.
Afterwards, we uniquely determine $M^{(a,b)}(t,\nu)$ via solving the ODE (\ref{phase:Neq})
 using a variant of the integral equation method of \cite{Greengard} 
 (or any standard method for stiff problems).

To obtain the
values of $\psi^{(a,b)}(t,\nu)$, we first calculate
the values of
\begin{equation}
\frac{d}{dt} \psi^{(a,b)}(t,\nu)
\end{equation}
via (\ref{phase:mpsip}).      Next, we obtain
the values of the function $\tilde{\psi}$  defined via
\begin{equation}
\tilde{\psi}(t) = \int_{\alpha_1}^t 
\frac{d}{ds} \psi^{(a,b)}(s,\nu)\ ds
\end{equation}
at any point in $[\frac{1}{N_{\mbox{\tiny max}}},\pi-\frac{1}{N_{\mbox{\tiny max}}}]$ via spectral integration.  There is
an unknown constant connecting $\tilde{\psi}$ with 
 $\psi^{(a,b)}(t,\nu)$; that is,
\begin{equation}
\psi^{(a,b)}(t,\nu) = \tilde{\psi}(t) + C.
\end{equation}
To evaluate $C$, we first use a combination of asymptotic and series
expansions to calculate $\tilde{P}_\nu^{(a,b)}$ at the point $\alpha_1$.
Since $\widetilde{\psi}(\alpha_1) =0$, it follows that 
\begin{equation}
\tilde{P}_\gamma^{(a,b)}(\alpha_1) = M^{(a,b)}(\alpha_1,\gamma) \cos( C),
\label{precomp:c}
\end{equation}
and $C$ can be calculated in the obvious fashion.

The above discussion has introduced an algorithm to evaluate $\psi^{(a,b)}(t,\nu)$ and $M^{(a,b)}(t,\nu)$ in the whole domain $(t,\nu)\in [\frac{1}{N_{\mbox{\tiny max}}},\pi-\frac{1}{N_{\mbox{\tiny max}}}]\times (27,N_{\mbox{\tiny max}}]$. To achieve a fast algorithm, we note that it is enough to conduct calculation on selected important grid points of $[\frac{1}{N_{\mbox{\tiny max}}},\pi-\frac{1}{N_{\mbox{\tiny max}}}]\times (27,N_{\mbox{\tiny max}}]$ through the above calculation, since we can evaluate $\psi^{(a,b)}(t,\nu)$ and $M^{(a,b)}(t,\nu)$ via interpolation with their values on the important grid points.

It is well known that a polynomial $p(t)$ (here $p(t)$ is either $\tilde{P}_\nu^{(a,b)}(t)$ or $\tilde{Q}_\nu^{(a,b)}(t)$) can be evaluated in a numerically stable fashion
using the barycentric Chebyshev interpolation formula \cite{Trefethen}
\begin{equation}
p(t) = \sum_{j=1}^{k}  \frac{w_j}{t-x_j} p(x_j)
\ \bigg / \ \ 
\sum_{j=1}^{k}  \frac{w_j}{t-x_j},
\label{cheby:bary}
\end{equation}
where $x_1,\ldots,x_k$ are the nodes of the $k$-point Chebyshev grid
on a sufficiently small interval $(\sigma_1,\sigma_2)$ (such that $p(t)$ is a polynomial of degree at most $k-1$) and 
\begin{equation}
w_j = \begin{cases}
(-1)^j,  & 1 < j < k; \\
(-1)^j\  \frac{1}{2}, & \mbox{otherwise}.
\end{cases}
\end{equation}
Hence, it is sufficient to evaluate $\psi^{(a,b)}(t,\nu)$ and $M^{(a,b)}(t,\nu)$, which lead to $\tilde{P}_\nu^{(a,b)}(t)$ and $\tilde{Q}_\nu^{(a,b)}(t)$, at a tensor product of Chebyshev grid points in $t$ and $\nu$, and calculate them at an arbitrary location $(t,\nu)$ via a bivariate piecewise barycentric Chebyshev interpolation.

More explicitly, 
let $m_\nu$ be the least integer such that $3^{m_\nu} \geq N_{\mbox{\tiny max}}$, and let $m_t$ be equal to twice the least integer $l$ such that $\frac{\pi}{2} 2^{-l+1} \leq \frac{1}{N_{\mbox{\tiny max}}}$. Next, we define  $\beta_j = \max\left\{ 3^{j+2}, N_{\mbox{\tiny max}}\right\}$ for $j=1,\ldots,m_\nu$, $\alpha_i = \frac{\pi}{2}  2^{i-m_t/2}$ for $i=1,\ldots,m_t/2$, and $\alpha_i= \pi - \frac{\pi}{2}  2^{m_t/2+1-i}$ for $i=m_t/2+1, \ldots,m_t$. Now we let 
\begin{equation}
\tau_1,\ldots,\tau_{M_t}
\label{precomp:tnodes}
\end{equation}
denote the $16$-point piecewise Chebyshev gird on the intervals
\begin{equation}
\left(\alpha_1,\alpha_2\right), \left(\alpha_2,\alpha_3\right), \ldots,
\left(\alpha_{m_t-1},\alpha_{m_t}\right)
\label{precomp:tints}
\end{equation}
with $M_t = 15( m_t-1) + 1$ points in total; let 
\begin{equation}
\gamma_1,\ldots,\gamma_{M_\nu}
\label{precomp:nunodes}
\end{equation}
denote the nodes of the $24$-point piecewise Chebyshev grid
on the intervals
\begin{equation}
\left(\beta_1,\beta_2\right), \left(\beta_2,\beta_3\right), \ldots,
\left(\beta_{{m_\nu}-1},\beta_{m_{\nu}}\right)
\label{precomp:nuints}
\end{equation}
with $M_t = 23( m_\nu-1) + 1$ points in total. The piecewise  Chebyshev grids \eqref{precomp:tnodes} and \eqref{precomp:nunodes} form a tensor product \begin{equation}
\label{eq:tsp}
\left\{
\left(\tau_i,\gamma_j\right) : i=1,\ldots,M_t, \ j=1,\ldots,M_\nu
\right\}.
\end{equation}
Then $\psi^{(a,b)}(t,\nu)$ and $M^{(a,b)}(t,\nu)$ can be obtained via barycentric Chebyshev interpolations in $t$ and $\nu$ via their values at the tensor product restricted to the piece $\left(\alpha_i,\alpha_{i+1}\right)\times \left(\beta_j,\beta_{j+1}\right)$ that $(t,\nu)$ belongs to for some $i$ and $j$.

Here we summarize the operation complexity of the whole algorithm above. A detailed discussion can be found in \cite{Jacobi}. ODE solvers are applied for 
$\mathcal{O}\left(\log\left(N_{\mbox{\tiny max}}\right)\right)$ 
values of $\nu$, and the solution of the solver is evaluated at
$\mathcal{O}\left(\log\left(N_{\mbox{\tiny max}}\right)\right)$ 
values of $t$, making the total running time
of the procedure just described $\mathcal{O}\left(\log^2\left(N_{\mbox{\tiny max}}\right)\right)$.

Once the values of $\psi_{\nu}^{(a,b)}$ and $M_\nu^{(a,b)}$
are ready at the tensor product
of the piecewise Chebyshev grids (\ref{precomp:tnodes}) and (\ref{precomp:nunodes}),
they can be evaluated for any $t$ and $\nu$ via repeated application
of the barycentric
Chebyshev interpolation formula  in
the same number of operations which
is independent of $\nu$ and $t$.

\section{Fast Jacobi polynomial transforms}

\label{sec:jt}

Without loss of generality, we assume that $\psi^{(a,b)}(t,\nu)$, $M^{(a,b)}(t,\nu)$, $P_\nu^{(a,b)}$, and $Q_\nu^{(a,b)}$ can be evaluated in $O(1)$ operations for any $t\in (0,\pi)$ and $\nu\in[0,N_{\mbox{\tiny max}} ]$. A fast algorithm for rapidly computing Gauss-Jacobi quadrature nodes and weights in the Jacobi polynomial transform in \eqref{eq:jactrs2} has also been introduced in \cite{Jacobi}. We refer the reader to \cite{Jacobi} for details and assume that quadrature nodes and weights are available in this section.

\subsection{One-dimensional transform and its inverse}
\label{sec:EOJT}
Here we propose a new variant of the one-dimensional Jacobi polynomial transform for $\tilde{P}_j^{(a,b)}(t)$ in \cite{Jacobi}. 
The new algorithm simplifies the discussion of the fast algorithm. The transform for $\tilde{Q}_j^{(a,b)}(t)$ is similar. 
For our purpose, we consider the $n^{th}$ order uniform forward Jacobi polynomial transform
 with a scaling, that is, calculate the vector  of values
\begin{equation}
\left(
\begin{array}{c}
f(t_1) \sqrt{w_1} \\
f(t_2) \sqrt{w_2} \\
\vdots\\
f(t_n) \sqrt{w_n} \\
\end{array}
\right)
\label{tf:transout}
\end{equation}
given the vector
\begin{equation}
\left(
\begin{array}{c}
\alpha_1 \\
\alpha_2\\
\vdots \\
\alpha_n
\end{array}
\right)
\label{tf:transin}
\end{equation}
of the coefficients in the expansion
\begin{equation}
f(t) = \sum_{j=0}^{n-1} \alpha_j \tilde{P}_j^{(a,b)}(t);
\label{tf:transexp}
\end{equation}
here,
 $t_1,\ldots,t_n,w_1\ldots,w_n$ are the nodes and weights
of the $n$-point trigonometric Gauss-Jacobi quadrature rule corresponding
to the parameters $a$ and $b$.  The ``nonuniform" forward Jacobi transform does not require trigonometric Gauss-Jacobi quadrature nodes and weights.

In the uniform transform, the properties of the trigonometric Gauss-Jacobi quadrature rule and the
weighting by square roots in (\ref{tf:transout}) ensure that the
$n \times n$ matrix $W \mathcal{J}_n^{(a,b)}$
taking (\ref{tf:transin}) to (\ref{tf:transout}) is orthogonal, where $W$ is the $n\times n$ matrix
\begin{equation}
W = \left(
\begin{array}{cccc}
\sqrt{w_1} & 0          & 0      & 0           \\
0          & \sqrt{w_2} & 0      & 0           \\
0          & 0          & \ddots & 0           \\
0          & 0          & 0      & \sqrt{w_n} \\
\end{array}
\right),
\end{equation}
and $\mathcal{J}_n^{(a,b)}$ is the $n \times n$ matrix representing the sum in \eqref{tf:transexp}. Hence, the inverse transform is a simple matvec by the transpose of $W\mathcal{J}_n^{(a,b)}$. In the nonuniform case, $\mathcal{J}_n^{(a,b)}$ is usually a very ill-conditioned matrix and the inverse transform requires solving a challenging system of linear equations, which will be discussed in a separate paper in the future.

It follows from (\ref{introduction:p}) 
 that the $(j,k)$-th entry of $\mathcal{J}_n^{(a,b)}$ 
is
\begin{align}
\label{tf:jn} 
 \left(\mathcal{J}_n^{(a,b)} \right)_{jk}=& \begin{cases}
   M^{(a,b)}(t_j,k-1) \cos\left( \psi^{(a,b)}(t_j,k-1) \right),
& \text{if } k-1>27, \\
   \tilde{P}_{k-1}^{(a,b)}  (t_j),      & \text{otherwise,}
  \end{cases}
\end{align}
since the part consisting of low-degree polynomials is computed via three-term recurrence relations or various asymptotic expansions. Let $\mathcal{V}_n^{(a,b)}\in\mathbb{R}^{n\times 27}$ be the submatrix of $\mathcal{J}_n^{(a,b)}$ corresponding to the part of polynomials with degree less than $27$, and $\mathcal{G}_n^{(a,b)}\in\mathbb{R}^{n\times n-27}$ be the rest of $\mathcal{J}_n^{(a,b)}$, then 
\[
\mathcal{J}_n^{(a,b)} = \left[ \mathcal{V}_n^{(a,b)},\mathcal{G}_n^{(a,b)} \right].
\] 
Note that $\mathcal{G}_n^{(a,b)} $ is the real part of a discrete Fourier integral transform that can be evaluated with a quasi-linear complexity by the algorithm of \cite{Candes-Demanet-Ying} for applying Fourier integral transforms and a special case of that in \cite{yang2018unified} for general oscillatory integral transforms. Motivated by this observation, the algorithm in \cite{Jacobi} constructs a low-rank matrix $\mathcal{A}_{n}^{(a,b)}\in\mathbb{C}^{n \times (n-27)}$ whose $(j,k)$-th  entry is 
\begin{equation}
 M^{(a,b)}(t_j, k+27)\exp\left( \ii\left(\psi^{(a,b)} (t_j, k+27) - (k+27) t_j\right) \right).
\label{introduction:matrix}
\end{equation}
Denote the  $n \times (n-27)$ nonuniform FFT matrix, whose $(j,k)$-th entry is 
\begin{equation}
\exp\left( i (k+27) t_j\right),
\label{introduction:nufft}
\end{equation}
as $\mathcal{N}_{n}$. 
Then 
\begin{equation}
\label{eq:JAN}
\mathcal{J}_{n}^{(a,b)} = \left[ \mathcal{V}_n^{(a,b)},\Re{(\mathcal{A}_{n}^{(a,b)}\otimes \mathcal{N}_{n})}\right],
\end{equation}
where ``$\otimes$" denotes the Hadamard product. 
It can be immediately seen from \eqref{eq:JAN} that $\mathcal{J}_{n}^{(a,b)}$ can be applied through a small number of 
NUFFTs \cite{Greengard,Townsend} and a simple direct summation for $ \mathcal{V}_n^{(a,b)}$ in a quasi-linear scaling of $n$.

More specifically, let a low-rank approximation of $\mathcal{A}_{n}^{(a,b)}$ be
\begin{equation}
\label{eq:approA}
\mathcal{A}_{n}^{(a,b)} \approx \sum\limits_{j=1}^{r}u_{j}v_{j}^{T},
\end{equation}
where $u_{j}$ and $v_{j}$ are column vectors for each $j$. 
Then the transform from (\ref{tf:transin}) to (\ref{tf:transout}) can be approximated by the following sum,
\begin{equation}
\label{eq:appromv}
W\mathcal{J}_{n}^{(a,b)}\alpha  \approx W\mathcal{V}_n^{(a,b)} \alpha_{\mathcal{V}} + W\Re(\sum\limits_{j=1}^{r}D_{u_{j}}\mathcal{N}_{n}D_{v_{j}}\alpha_{\mathcal{G}}),
\end{equation}
where $\alpha_{\mathcal{V}}$ is the subvector of $\alpha$ with the first $27$ entries, $\alpha_{\mathcal{G}}$ is the subvector of $\alpha$ with the rest entries, $D_{u}$ denotes a diagonal matrix with a column vector $u$ on its diagonal. The formula (\ref{eq:appromv}) 
 can be carried out by $r$ NUFFTs in $O(rn\log n)$ arithmetic operations. {{Generally speaking, $r$ is related to $n$ and it was conjectured in \cite{Jacobi} that 
\begin{equation*}
r = O(\dfrac{\log n}{\log \log n})
\end{equation*} 
via excessive numerical experiments}}.

In practice, inspired by the NUFFT in \cite{Townsend}, we can replace the Hadamard product of a low-rank matrix and a NUFFT matrix in \eqref{eq:JAN} with a Hadamard product of a low-rank matrix and a DFT matrix, the matvec of which can be carried out more efficiently with a few numbers of FFTs.

Consider another matrix $\mathcal{B}_{n}^{(a,b)}\in\mathbb{C}^{n \times (n-27)}$ whose $(j,k)$-th entry is defined via 
\begin{equation}
\label{eq:B}
M^{(a,b)}\left(t_{j},k+27\right)\exp\left(\ii(\psi^{(a,b)}(t_{j},k+27)-2\pi\dfrac{[t_{j}n/2\pi]}{n}(k+27))\right),
\end{equation}
where $[x]$ denotes the integer nearest to $x$. Then we have
\begin{equation}
\label{eq:lastJ}
\mathcal{J}_{n}^{(a,b)} =  \left[ \mathcal{V}_n^{(a,b)},\Re{(\mathcal{B}_{n}^{(a,b)}\otimes \mathcal{F}_{n})}\right],
\end{equation}
where $\mathcal{F}_{n}$ is an $n\times (n-27)$ matrix whose $(j,k)$-th entry is
\begin{equation}
\label{eq:permuDFT}
\exp\left(i(2\pi\dfrac{[t_{j}n/2\pi]}{n}(k+27))\right).
\end{equation}
It's obvious that $\mathcal{F}_{n}$ in (\ref{eq:lastJ}) 
is a row permutation of an inverse DFT matrix. 
Note that the difference between the phase functions of $\mathcal{A}_{n}^{(a,b)}$ and $\mathcal{B}_{n}^{(a,b)}$ is less than $\pi$. Hence, $\mathcal{B}_{n}^{(a,b)}$ is also a low-rank matrix and we can apply the randomized SVD in Algorithm \ref{algo} to construct a low-rank approximation of $\mathcal{B}_{n}^{(a,b)}$ in $O(r^2n)$ operations.

Suppose we have constructed an approximate rank-$r$ SVD of $\mathcal{B}_{n}^{(a,b)}$ up to a desired accuracy using Algorithm \ref{algo}; that is, 
suppose that we have computed the factorization
\begin{equation}
\mathcal{B}_{n}^{(a,b)} \approx USV,
\end{equation}
where $U\in \bbC^{n\times r}$, $V\in \bbC^{r\times (n-27)}$, and $S\in \bbR^{r\times r}$ is a positive definite diagonal matrix. By rearranging the factors above, we have
\begin{equation}
\label{eq:approxA}
\mathcal{B}_{n}^{(a,b)} \approx (US^{\frac{1}{2}})(S^{\frac{1}{2}}V) = u_{1}v_{1}^{T}+\cdots+u_{r}v_{r}^{T},
\end{equation} 
where $u_{i}$ and $v_{i}^{T}$ denote the $i^{th}$ column vector of $US^{\frac{1}{2}}$ and the $i^{th}$ row vector of $S^{\frac{1}{2}}V$, respectively, and $T$ denotes the matrix transpose. Once (\ref{eq:approxA}) is ready, we have
\begin{equation}
\label{eq:oneJ}
W\mathcal{J}_{n}^{(a,b)}\alpha \approx W\mathcal{V}_n^{(a,b)} \alpha_{\mathcal{V}} + W\Re(((\sum\limits_{j=1}^{r}u_{j}v_{j}^T)\otimes \mathcal{F}_{n}) \alpha_{\mathcal{G}}) =W\mathcal{V}_n^{(a,b)} \alpha_{\mathcal{V}} + W\Re(\sum\limits_{j=1}^{r}D_{u_{j}}\mathcal{F}_{n}D_{v_{j}}\alpha_{\mathcal{G}}),
\end{equation}
where $D_{u}$ denotes a diagonal matrix with $u$ on its diagonal. Formula (\ref{eq:oneJ}) indicates that the Jacobi polynomial transform can be evaluated efficiently via $r$ inverse FFTs, which requires $O(rn\log n)$ arithmetic operations and $O(rn)$ memory. {{SVD gives an optimal numerical rank of a matrix given a fixed accuracy in the low-rank approximation. Hence, the rank in our new method should be smaller than or equal to the rank in \cite{Jacobi}. Numerical performance in this paper and \cite{Jacobi} leads to the conjecture that the rank 
\begin{equation*}
r = O(\dfrac{\log n}{\log \log n}).
\end{equation*} }} 

 Compared to (\ref{eq:appromv}) used in the original fast Jacobi transform in \cite{Jacobi}, the number of inverse FFTs in (\ref{eq:oneJ}) has been optimized. Thus Formula (\ref{eq:oneJ}) would take fewer operations to compute a multi-dimensional transform. Note that, compared to the original method in \cite{Jacobi}, though our new method using \eqref{eq:oneJ} would cost more time to construct a low-rank approximation, the optimized rank of the low-rank approximation would accelerate the application of the multi-dimensional Jacobi transform. Besides, as we mentioned previously, the new method works in a larger range of $a$ and $b$. In Section \ref{sec:results}, we will provide numerical comparisons to demonstrate the superiority of the new method over the method in \cite{Jacobi}.

The fast inverse Jacobi transform in the uniform case can be carried out in a similar manner, since $W\mathcal{J}_{n}^{(a,b)}$ is an orthonormal matrix. In fact, the inverse transform can be computed via
\[
\alpha_{\mathcal{V}} = \left(\mathcal{V}_n^{(a,b)}\right)^T W^TWf,
\]
and
\begin{equation}
\label{eq:oneJinv}
\alpha_{\mathcal{G}}\approx \Re(((\sum\limits_{j=1}^{r}u_{j}v_{j}^T)\otimes \mathcal{F}_{n})^{T} W^{T}Wf) = \Re(\sum\limits_{j=1}^{r}D_{v_{j}}\mathcal{F}_{n}^{T}D_{u_{j}}W^{T}Wf),
\end{equation}
where $\mathcal{F}_{n}^{T}$ is a permutation of the DFT matrix. Therefore, the inverse Jacobi polynomial transform can also be computed via $r$ FFTs.

It is worth emphasizing that the fast algorithm introduced in this section also works for non-uniform Jacobi polynomial transforms since the low-rankness of $\mathcal{B}_{n}^{(a,b)}$ is independent of the samples in $t$ and the quadrature weights in the uniform Jacobi polynomial transform. There will be numerical examples later to verify this. The inverse transform in the non-uniform case requires solving a highly ill-conditioned linear system, which will be reported in a separate paper in the future.

\subsection{Two-dimensional transform and its inverse} 
\label{sec:DMJT}
In this section, we extend our algorithm to the two-dimensional case, using the new method developed in the last section. We only focus on the transform for $\tilde{P}_{\nu}^{(a,b)}(t)$; the one for $\tilde{Q}_{k-1}^{(a,b)}(t)$ is similar. {{Furthermore, we concentrate on transforms admitting a tensor-product structure which is usually the case in real  applications.}} We will adopt the MATLAB notation $A(:)$ or $vec(A)$ as a vector resulting from reshaping the matrix $A$ into a vector. We also use similar notations as in Section \ref{sec:EOJT}, e.g., $\mathcal{J}_{n,\cdot}^{(a,b)}$, $\mathcal{V}_{n,\cdot}^{(a,b)}$, $\mathcal{G}_{n,\cdot}^{(a,b)}$, $\mathcal{B}_{n,\cdot}^{(a,b)}$, $\mathcal{F}_{n,\cdot}$, and $W_{\cdot}$ denote corresponding matrices analogous to their counterparts in Section \ref{sec:EOJT}, respectively, with ``$\cdot$" specifying a variable $x$ or $y$  in the spatial domain. In the rest of this section, we always assume a low-rank approximation 
\begin{equation}
\label{eqn:Blr}
\mathcal{B}_{n,\cdot}^{(a,b)} = \sum_{i=1}^{r_{\cdot}}u_{\cdot,i}v_{\cdot,i}^{T}
\end{equation}
has been obtained by Algorithm \ref{algo} up to a desired accuracy for ``$\cdot$" as $x$ or $y$.

Given locations $\{x_{i}\}_{i=1,\cdots,n} \subset (0,\pi)$ and $\{y_{i}\}_{i=1,\cdots,n} \subset (0,\pi)$, with no substantial difference, the forward and inverse two-dimensional Jacobi
polynomial transforms arise in the following Jacobi expansion
\begin{equation}
\label{eq:2Djac}
f(x_i,y_j) = \sum\limits_{k=1}^{n}\sum\limits_{\ell=1}^{n}\alpha(k,\ell)\tilde{P}_{k-1}^{(a,b)}(x_{i})\tilde{P}_{\ell-1}^{(a,b)}(y_{j}),\quad\text{for } i,j = 1,\cdots,n,
\end{equation}
where $\alpha$ denotes an expansion coefficients matrix. The forward and inverse transform can be defined analogously as in the one-dimensional case. When both $\{x_{i}\}_{i=1,\cdots,n}$ and $\{y_{i}\}_{i=1,\cdots,n}$ are exactly the nodes of the trigonometric Gauss-Jacobi quadrature rule, the corresponding transform is referred to as the uniform transform. With additional weight matrices $W_x$ and $W_y$, the forward transformation matrix $(W_{x}\mathcal{J}_{n,x}^{(a,b)})\odot (W_{y}\mathcal{J}_{n,y}^{(a,b)})$ taking $\alpha(:)$ to $(W_{x}\odot W_{y})f(:)$ is orthogonal, where ``$\odot$" denotes the Kronecker product. When $\{x_{i}\}_{i=1,\cdots,n}$, $\{y_{i}\}_{i=1,\cdots,n}$, and weights are not given by the trigonometric Gauss-Jacobi quadrature rule, we refer the corresponding transform as a non-uniform transform.

With the low-rank approximations in \eqref{eqn:Blr} available, {{we have

\begin{equation}
\label{eq:TJ}
\begin{split}
\left[(W_{x}\mathcal{J}_{n,x}^{(a,b)})\odot (W_{y}\mathcal{J}_{n,y}^{(a,b)})\right]\alpha(:) = vec \left[((W_{x}\mathcal{J}_{n,x}^{(a,b)})((W_{y}\mathcal{J}_{n,y}^{(a,b)})\alpha)^{T})^{T} \right].
\end{split}
\end{equation}

This formula indicates that 2D uniform Jacobi polynomial transforms can be evaluated via $2n$ 1D forward transforms, thus $O((r_{x}+r_{y})n)$ 1D inverse FFTs, which results in $O((r_{x}+r_{y})n^2 \log(n))$ arithmetic operations and $O((r_{x}+r_{y})n^2)$ memory. We would like to emphasize that this algorithm also works for a 2D non-uniform forward transform since it is a tensor form of the one-dimensional transform. There will be numerical
examples later to verify this. Note that although the left-hand side and the right-hand side of Formula (\ref{eq:TJ}) are mathematically equivalent, we adopt the right-hand side in implementation to compute the 2D transform because Level 3 BLAS can be applied to accelerate matrix-matrix multiplications. Analogous techniques can be used to compute other multidimensional transforms in this paper. }}

The fast 2D inverse Jacobi transform in the uniform case can be carried out in a similar manner, since both $W_{x}\mathcal{J}_{n,x}^{(a,b)}$ and $W_{y}\mathcal{J}_{n,y}^{(a,b)}$ are orthonormal matrices. {{In fact, the 2D inverse transform can be computed via

\begin{equation}
\label{eq:invTJ}
\alpha(:) = \left[(\mathcal{J}_{n,x}^{(a,b)})^T\odot (\mathcal{J}_{n,y}^{(a,b)})^{T}\right]\tilde{f}(:) = vec((\mathcal{J}_{n,y}^{(a,b)})^{T}\tilde{f}\mathcal{J}_{n,x}^{(a,b)}),
\end{equation}
where $\tilde{f} := W_{y}^{2} f W_{x}^{2}$.
Again, Level 3 BLAS and $2n$ 1D inverse transforms can be applied to compute this 2D transform.}} Therefore, 2D inverse Jacobi polynomial transforms can also be evaluated via $O((r_{x}+r_{y})n)$ 1D FFTs.

\subsection{Three-dimensional transform and its inverse} 
\label{sec:TMJT}
In this section, we continue to extend our algorithm for 3D Jacobi polynomial transform and its inverse. Analogously, we just discuss the transform for $\tilde{P}_{\nu}^{(a,b)}(t)$ and the transform for $\tilde{Q}_{k-1}^{(a,b)}(t)$ is similar. {{Just like the 2D case, we mainly focus on those transforms with a tensor-product structure.}} The notations in Section \ref{sec:DMJT} will be inherited here.

Given locations $\{x_{i}\}_{i=1,\cdots,n}$, $\{y_{i}\}_{i=1,\cdots,n}$, and $\{z_{i}\}_{i=1,\cdots,n}$ in $(0,\pi)$, with no substantial difference, the forward and inverse three-dimensional Jacobi
polynomial transforms arise in the following Jacobi expansion
\begin{equation}
\label{eq:3Djac}
f(x_i,y_j,z_k) = \sum\limits_{h=1}^{n}\sum\limits_{\ell=1}^{n}\sum\limits_{m=1}^{n}\alpha(k,\ell,m)\tilde{P}_{h-1}^{(a,b)}(x_{i})\tilde{P}_{\ell-1}^{(a,b)}(y_{j})\tilde{P}_{m-1}^{(a,b)}(z_{k}),\quad\text{for } i,j,k = 1,\cdots,n,
\end{equation}
where $\alpha$ denotes a three-dimensional tensor containing expansion coefficients. The forward and inverse transforms can be defined analogously to the three-dimensional case. When $\{x_{i}\}_{i=1,\cdots,n}$, $\{y_{i}\}_{i=1,\cdots,n}$, and $\{z_{i}\}_{i=1,\cdots,n}$ are all exactly the nodes of the trigonometric Gauss-Jacobi quadrature rule, the corresponding transform is referred to as the uniform transform. Otherwise, we refer the corresponding transform as a non-uniform transform.

In the uniform transform, in order to take advantage of orthogonality, we consider the tensor product $(W_{x}\mathcal{J}_{n,x}^{(a,b)})\odot (W_{y}\mathcal{J}_{n,y}^{(a,b)})\odot (W_{y}\mathcal{J}_{n,z}^{(a,b)})$ that takes $\alpha(:)$ to $(W_{x}\odot W_{y}\odot W_{z})f(:)$. Once the related low-rank approximations, 
\begin{equation}
\label{eqn:Blr2}
\mathcal{B}_{n,\cdot}^{(a,b)} = \sum_{i=1}^{r_{\cdot}}u_{\cdot,i}v_{\cdot,i}^{T},
\end{equation}
have been obtained, we have
\begin{equation}
\label{eq:TJ3pre}
(W_{x}\odot W_{y}\odot W_{z})f(:) = (W_{x}\odot W_{y}\odot W_{z})(\mathcal{J}_{n,x}^{(a,b)}\odot \mathcal{J}_{n,y}^{(a,b)}\odot \mathcal{J}_{n,z}^{(a,b)})\alpha(:),
\end{equation}
{{where

\begin{equation}
\label{eq:TJ3}
(\mathcal{J}_{n,x}^{(a,b)}\odot \mathcal{J}_{n,y}^{(a,b)}\odot \mathcal{J}_{n,z}^{(a,b)})\alpha(:) = vec\left[((\mathcal{J}_{n,y}^{(a,b)}\odot \mathcal{J}_{n,z}^{(a,b)})\tilde{\alpha}) (\mathcal{J}_{n,x}^{(a,b)})^{T}\right]
\end{equation}
and $\tilde{\alpha}$  is a $n^2 \times n$ matrix reshaped from $\alpha$.

Formula (\ref{eq:TJ3}) implies that a 3D uniform Jacobi polynomial transform can be computed by evaluating the right-hand side using Level 3 BLAS routines with $n$ 2D forward transforms and $n^2$ 1D forward transforms, which are dominated by totally $O((r_{x}+r_{y}+r_{z})n^2)$ inverse FFTs and result in $O((r_{x}+r_{y}+r_{z})n^3 \log(n))$ arithmetic operations and $O((r_{x}+r_{y}+r_{z})n^3)$ memory. Again, we would like to emphasize that this algorithm also works for 3D non-uniform forward transforms.}} There will be numerical
examples later to verify this.

The fast 3D inverse (uniform) Jacobi transform in the uniform case can be carried out in a similar manner, since $W_{x}\mathcal{J}_{n,x}^{(a,b)}$, $W_{y}\mathcal{J}_{n,y}^{(a,b)}$ and $W_{z}\mathcal{J}_{n,z}^{(a,b)}$ are all orthonormal matrices. In fact, the 3D inverse transform can be computed via
{{
\begin{equation}
\label{eq:invTJ3}
(\mathcal{J}_{n,x}^{(a,b)}\odot \mathcal{J}_{n,y}^{(a,b)}\odot \mathcal{J}_{n,z}^{(a,b)})^{T}\tilde{f}(:) = vec\left[(((\mathcal{J}_{n,y}^{(a,b)})^{T}\odot (\mathcal{J}_{n,z}^{(a,b)})^{T})\bar{f}) \mathcal{J}_{n,x}^{(a,b)}\right]
\end{equation}
where $\tilde{f}(:) = (W_{x}\odot W_{y}\odot W_{z})^2 f(:)$ and $\bar{f}$ is an $n^2\times n$ matrix reshaped from $\tilde{f}$.

Analogously, the right-hand side allows Level 3 BLAS with $n$ inverse 2D transforms and $n^2$ 1D inverse transforms.  Therefore, the 3D inverse Jacobi polynomial transform can also be evaluated via $O((r_{x}+r_{y}+r_{z})n^2)$ 1D FFTs, resulting in a nearly linear scaling algorithm.}}


\begin{figure}
  \begin{center}
    \begin{tabular}{cc}
      \includegraphics[height=1.45in]{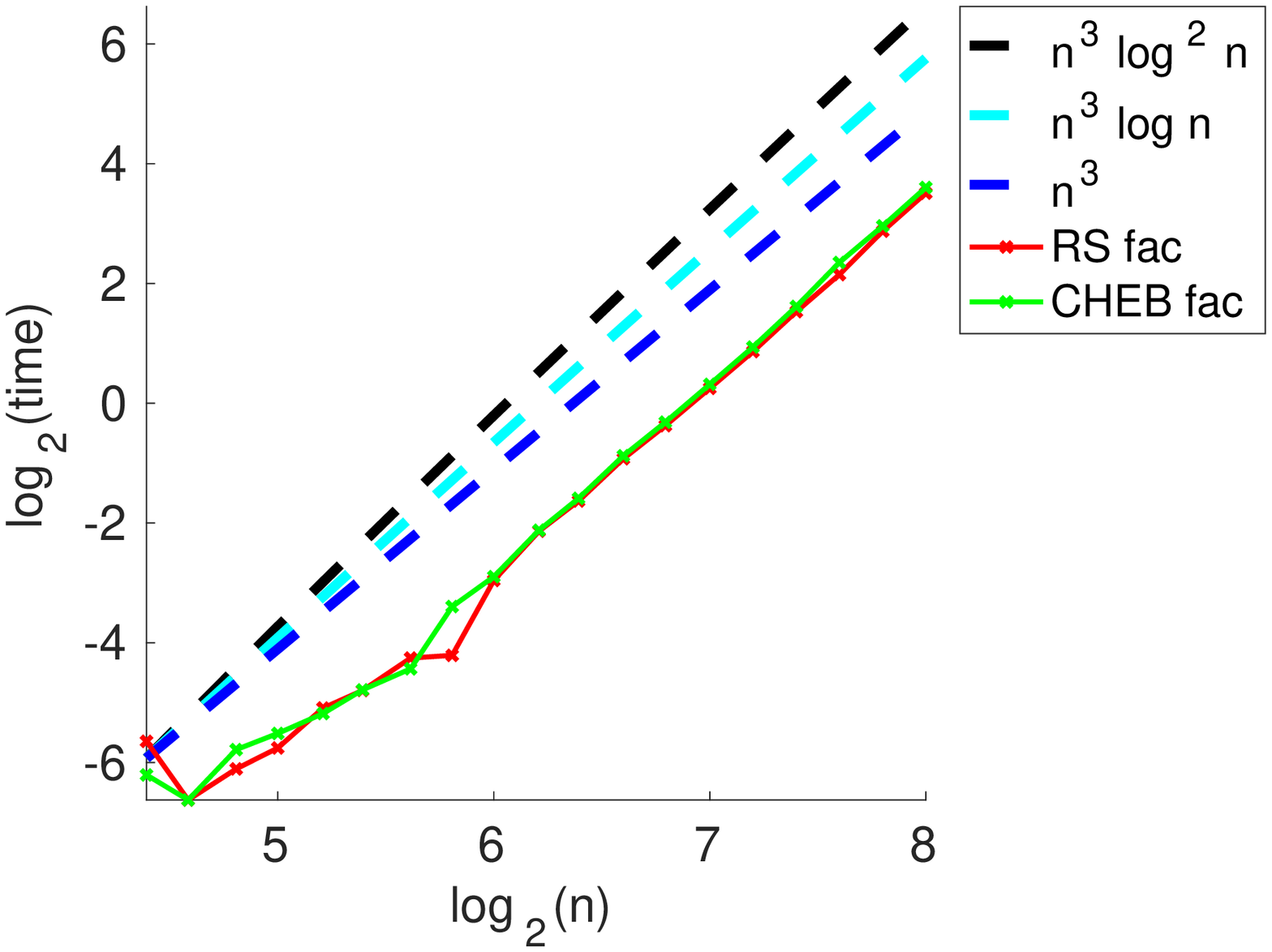}&
      \includegraphics[height=1.45in]{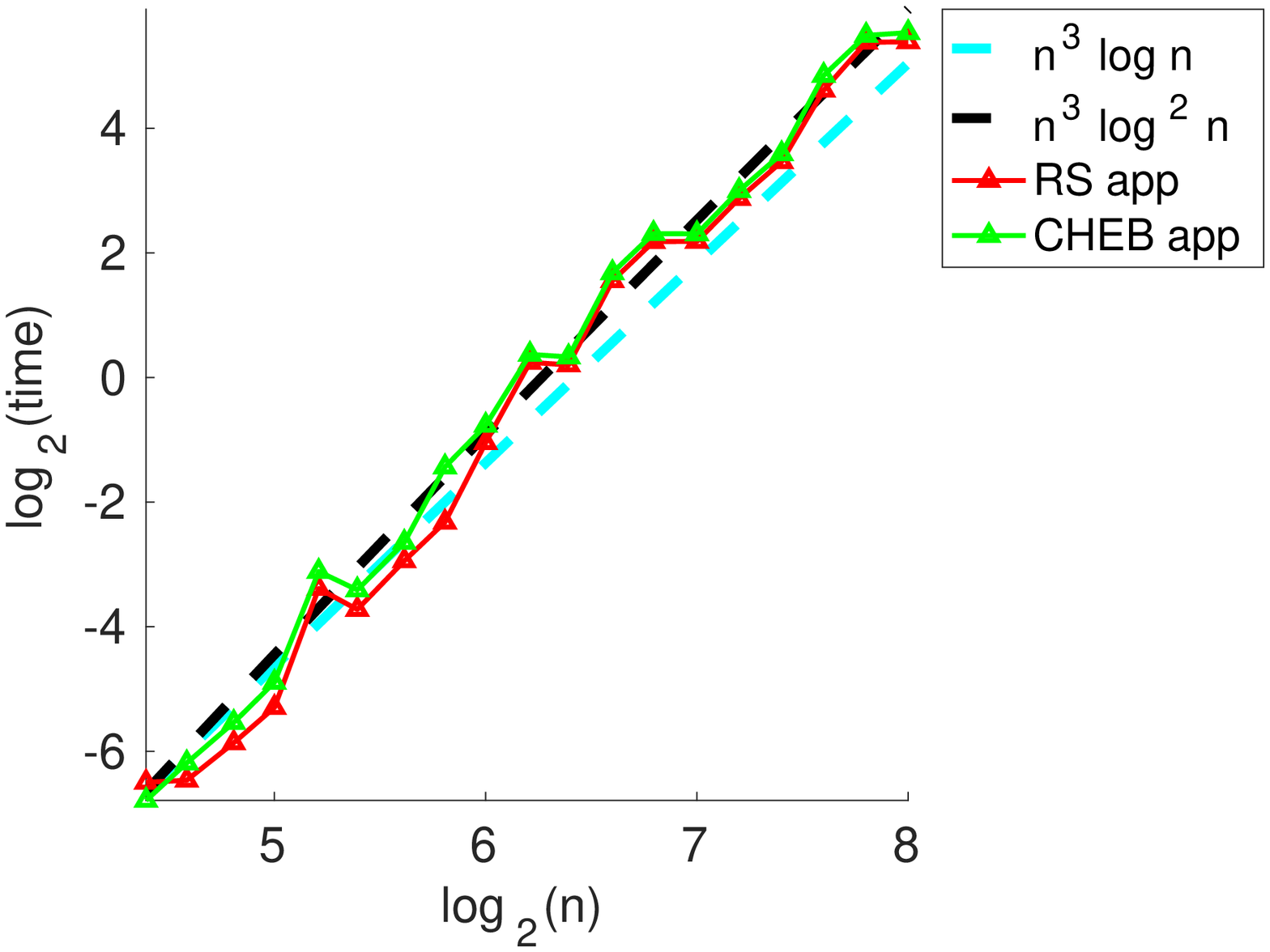}\\
      \\
      \includegraphics[height=1.45in]{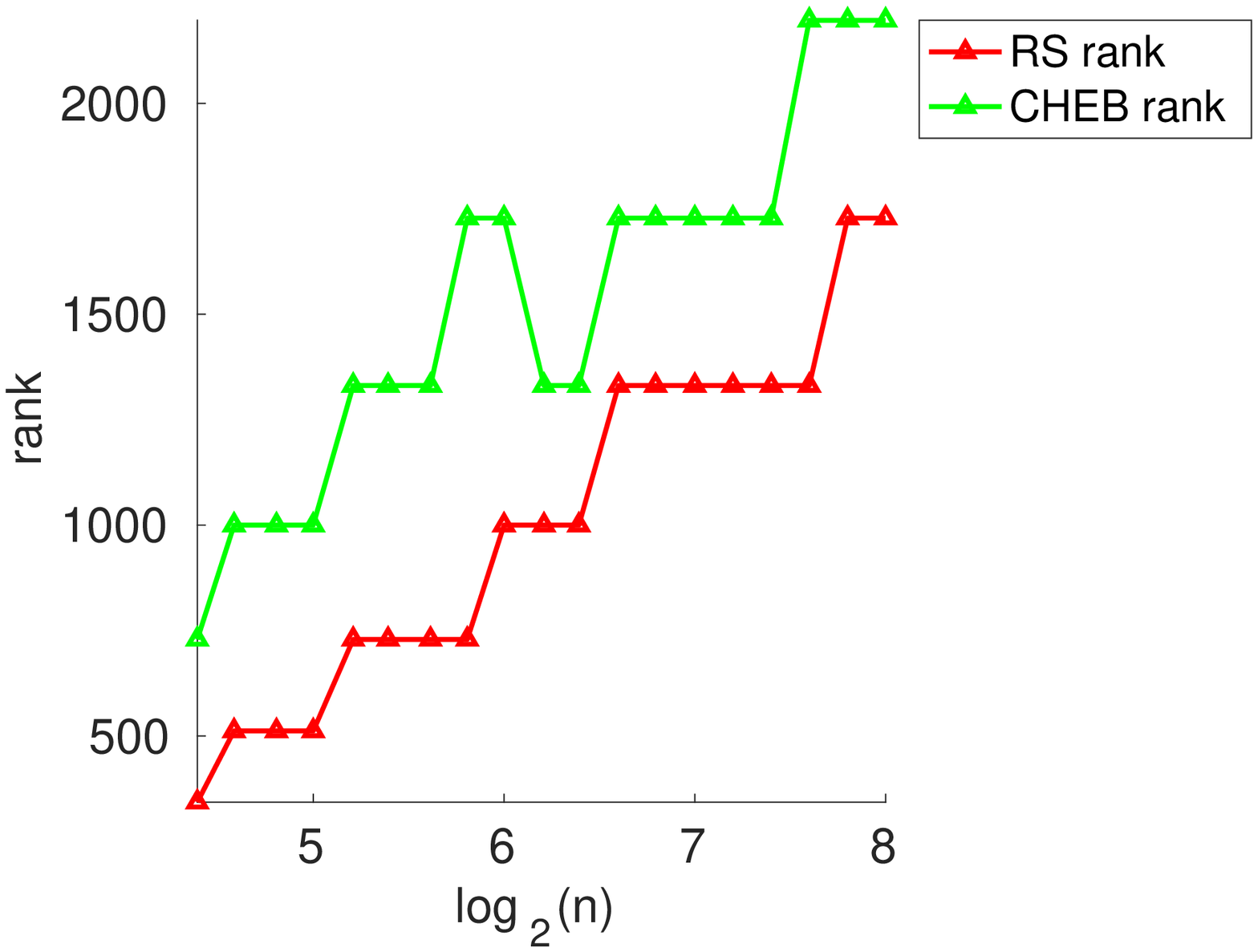}&
      \includegraphics[height=1.45in]{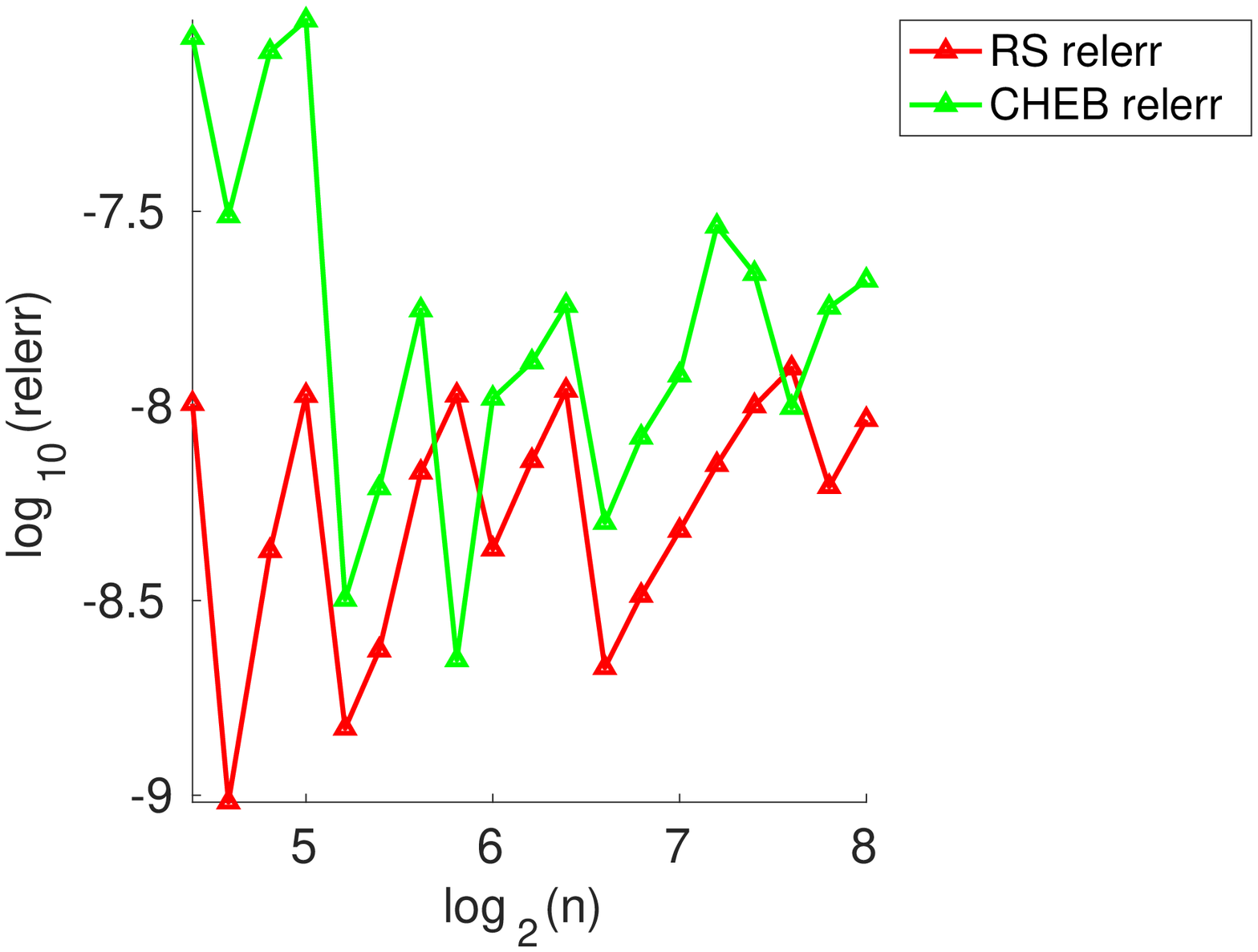}\\
    \end{tabular}
  \end{center}
\caption{Numerical results for the comparison of the CHEB and RS method for 3D forward uniform transform. Top: the factorization and application time {{(in second)}} from left to right, respectively. Bottom left and right: the numerical ranks of the low-rank matrices provided by CHEB and RS, and the relative errors of the fast matvec compared to the exact summation.}
\label{fig:comp3}
\end{figure}


\begin{table}
\begin{center}
\begin{tabular}{c|ccc|ccc|ccc}
\hline
\multirow{2}{*}{\diagbox{a=b}{$\log_{2}n$}} & \multicolumn{3}{c|}{1D} & \multicolumn{3}{c|}{2D} & \multicolumn{3}{c}{3D} \\ \cline{2-10}
& 10 & 15 & 20 & 6 & 9 & 12 & 6 & 7 & 8 
\\
\hline
-0.75 & 1.00 & 13.7 & 42.3 & 2.50 & 0.10 & 0.005 & 3.05 & 0.89 & 2.36 \\
\hline
-0.50 & 0.001 & 0.23 & 4.25 & 1E-5 & 1E-4 & 0.001 & 1E-4 & 7E-5 & 7E-5 \\
\hline
-0.25 & 0.33 & 4.65 & 44.3 & 1.57 & 0.07 & 0.005 & 1.87 & 1.04 & 3.56 \\
\hline
0.00 & 0.69 & 8.10 & 60.0 & 3.16 & 0.24 & 0.009 & 3.81 & 0.99 & 1.22 \\
\hline
0.25 & 0.71 & 1.95 & 23.7 & 4.74 & 0.36 & 0.009 & 5.90 & 1.52 & 1.66 \\
\hline
0.50 & 0.30 & 0.34 & 4.34 & 0.03 & 1E-4 & 0.001 & 0.97 & 0.23 & 0.35 \\
\hline
0.75 & 2.06 & 4.80 & 60.1 & 5.81 & 0.37 & 0.02 & 7.13 & 2.13 & 1.41 \\
\hline

\end{tabular}
\caption{Numerical results of Equation (\ref{eq:stability}) in the cases of 1D, 2D and 3D to measure stability. All numbers in this table are measured in the unit of $1E-8$. Hence, the quantity in Equation (\ref{eq:stability})  is of order $1E-8$ in all test cases. }
\label{tb:stability}
\end{center}
\end{table}

\begin{table}
\begin{center}
\begin{tabular}{c|c|c|c|c|c|c}
\hline
\diagbox{a=b}{rank}{$\log_{2}n$} & 14 & 15 & 16 & 17 & 18 & 19 \\
\hline
0.0 & $18,20$ & $19,21$ & $19,22$ & $20,23$ & $21,24$ & $20,26$ \\
\hline
0.1 & $19,20$ & $19,21$ & $19,23$ & $20,23$ & $21,24$ & $20,26$ \\
\hline
0.2 & $18,20$ & $18,21$ & $19,22$ & $20,23$ & $21,24$ & $20,26$ \\
\hline
0.3 & $18,19$ & $19,21$ & $20,22$ & $19,23$ & $21,24$ & $20,25$ \\
\hline
0.4 & $18,19$ & $19,19$ & $19,21$ & $19,22$ & $21,24$ & $19,25$ \\
\hline
0.5 & $9,9$ & $9,9$ & $9,9$ & $9,9$ & $9,9$ & $9,9$ \\
\hline
0.6 & $18,19$ & $18,20$ & $19,21$ & $19,23$ & $20,23$ & $20,24$ \\
\hline
0.7 & $18,19$ & $18,20$ & $20,21$ & $20,22$ & $21,23$ & $20,24$ \\
\hline
0.8 & $18,19$ & $19,20$ & $20,21$ & $19,21$ & $20,23$ & $20,24$ \\
\hline
0.9 & $18,18$ & $19,20$ & $20,21$ & $19,21$ & $20,22$ & $19,23$ \\
\hline

-0.1 & $18,20$ & $19,21$ & $19,23$ & $19,23$ & $20,25$ & $20,26$ \\
\hline
-0.2 & $17,20$ & $18,21$ & $18,22$ & $19,23$ & $20,24$ & $18,27$ \\
\hline
-0.3 & $18,19$ & $18,21$ & $18,22$ & $18,23$ & $20,24$ & $18,25$ \\
\hline
-0.4 & $17,19$ & $17,19$ & $17,21$ & $18,22$ & $19,24$ & $18,25$ \\
\hline
-0.5 & $2,9$ & $2,9$ & $3,9$ & $3,9$ & $4,9$ & $4,9$ \\
\hline
-0.6 & $16,19$ & $18,20$ & $17,21$ & $18,22$ & $18,23$ & $18,24$ \\
\hline
-0.7 & $17,19$ & $18,20$ & $18,21$ & $18,22$ & $19,23$ & $18,24$ \\
\hline
-0.8 & $17,21$ & $18,21$ & $18,22$ & $18,23$ & $19,25$ & $19,25$ \\
\hline
-0.9 & $18,20$ & $18,21$ & $18,22$ & $18,22$ & $19,23$ & $19,24$ \\
\hline

\end{tabular}
\caption{Numerical ranks of the CHEB and RS algorithms for 1D forward uniform transform up to a fixed accuracy $1E-8$ when $a=b$. There are two numbers in each item: the left one denotes the numerical rank by the RS algorithm and the right one is for the CHEB method in \cite{Jacobi}. }
\label{table1}
\end{center}
\end{table}

\begin{table}
\begin{center}
\begin{tabular}{c|c|c|c|c|c|c}
\hline
\diagbox{a=b}{error}{$\log_{2}n$} & 14 & 15 & 16 & 17 & 18 & 19 \\
\hline
0.0 & $0.17,0.66$ & $0.28,0.41$ & $0.13,0.77$ & $0.13,0.67$ & $0.95,0.57$ & $0.80,0.28$ \\
\hline
0.1 & $0.10,0.60$ & $0.45,0.46$ & $1.58,0.20$ & $0.73,0.51$ & $1.04,1.06$ & $5.47,0.45$ \\
\hline
0.2 & $0.48,0.25$ & $1.01,0.73$ & $1.51,0.77$ & $0.96,0.41$ & $0.73,1.08$ & $8.97,0.22$ \\
\hline
0.3 & $0.91,0.69$ & $0.32,0.47$ & $0.37,0.25$ & $0.15,0.43$ & $0.92,0.31$ & $2.41,0.80$ \\
\hline
0.4 & $0.30,1.19$ & $0.34,1.99$ & $0.71,1.21$ & $0.62,1.57$ & $0.52,0.53$ & $4.38,0.24$ \\
\hline
0.5 & $0.31,0.32$ & $0.31,0.27$ & $0.17,0.35$ & $0.11,0.14$ & $0.32,0.45$ & $0.17,0.27$ \\
\hline
0.6 & $0.32,0.40$ & $0.92,0.41$ & $1.51,0.33$ & $0.94,0.18$ & $3.98,0.76$ & $2.53,1.72$ \\
\hline
0.7 & $0.57,0.40$ & $1.03,1.40$ & $0.61,1.94$ & $1.72,1.17$ & $2.42,1.58$ & $4.07,2.30$ \\
\hline
0.8 & $0.79,0.93$ & $0.73,0.53$ & $0.83,1.43$ & $4.85,3.09$ & $7.15,1.45$ & $12.5,1.58$ \\
\hline
0.9 & $1.03,2.92$ & $1.01,0.84$ & $0.80,2.36$ & $3.57,4.84$ & $4.40,2.62$ & $26.9,5.43$ \\
\hline

-0.1 & $0.37,0.58$ & $0.69,0.33$ & $1.09,0.22$ & $1.80,0.59$ & $1.46,0.37$ & $7.44,0.39$ \\
\hline
-0.2 & $0.45,0.40$ & $0.33,0.41$ & $1.53,0.97$ & $1.65,0.48$ & $1.39,0.81$ & $9.30,0.21$ \\
\hline
-0.3 & $0.08,0.72$ & $0.36,0.33$ & $0.87,0.18$ & $1.54,0.40$ & $1.88,0.17$ & $16.4,0.95$ \\
\hline
-0.4 & $0.12,0.98$ & $0.64,0.93$ & $1.57,0.85$ & $2.03,0.97$ & $1.27,0.30$ & $4.22,0.14$ \\
\hline
-0.5 & $0.03,0.47$ & $0.05,0.40$ & $0.05,0.51$ & $0.31,0.31$ & $0.03,0.48$ & $0.08,0.47$ \\
\hline
-0.6 & $0.57,0.15$ & $0.17,0.47$ & $1.44,0.15$ & $0.99,0.15$ & $1.33,0.27$ & $3.97,1.79$ \\
\hline
-0.7 & $0.40,0.19$ & $0.45,1.72$ & $1.61,1.18$ & $3.13,1.68$ & $1.56,0.83$ & $5.56,3.62$ \\
\hline
-0.8 & $0.68,3E06$ & $1.43,1E06$ & $2.55,8E05$ & $3.34,7E05$ & $4.19,5E05$ & $4.87,3E05$ \\
\hline
-0.9 & $0.62,3E06$ & $1.25,1E06$ & $2.97,8E05$ & $3.21,6E05$ & $5.31,5E05$ & $8.73,3E05$ \\
\hline

\end{tabular}
\caption{Numerical relative accuracy of the CHEB and RS algorithms for 1D forward uniform transform up to a fixed accuracy $1E-8$ when $a=b$. The units of all data in this table are $1E-8$. There are two floats in each item: the left one denotes the relative accuracy of the low-rank factorization by the RS algorithm and the right one is for the CHEB algorithm in \cite{Jacobi}. }
\label{table2}
\end{center}
\end{table}

\section{Numerical results}
\label{sec:results}
This section presents several numerical examples to demonstrate the effectiveness of the algorithms
proposed above. Seaction \ref{sec:STA} demonstrates the stability of our algorithm. {Section \ref{sec:CDLFA} provides a comparison of the original method in \cite{Jacobi} named as method \textbf{CHEB} using \eqref{eq:appromv} in Section \ref{sec:EOJT}, our new method named as \textbf{RS} using \eqref{eq:oneJ} in Section \ref{sec:EOJT} to demonstrate the superiority of our algorithm. In Section \ref{sec:PFMT}, we apply our new method to 2D and 3D Jacobi polynomial transforms with parameters $a=b=0.40$ to show the complexity of our algorithm.  All implementations are in MATLAB\textsuperscript{\textregistered} on a server computer with 28 processors and 2.6 GHz CPU. {{But only one processor was used in every single experiment.  We have made our code, including that for all of the experiments described here,
available on GitLab at the following address:
  \begin{center}
\url{https://gitlab.com/FastOrthPolynomial/Jacobi.git}
  \end{center}

    \begin{figure}
  \begin{center}
    \begin{tabular}{ccc}
      \includegraphics[height=1.45in]{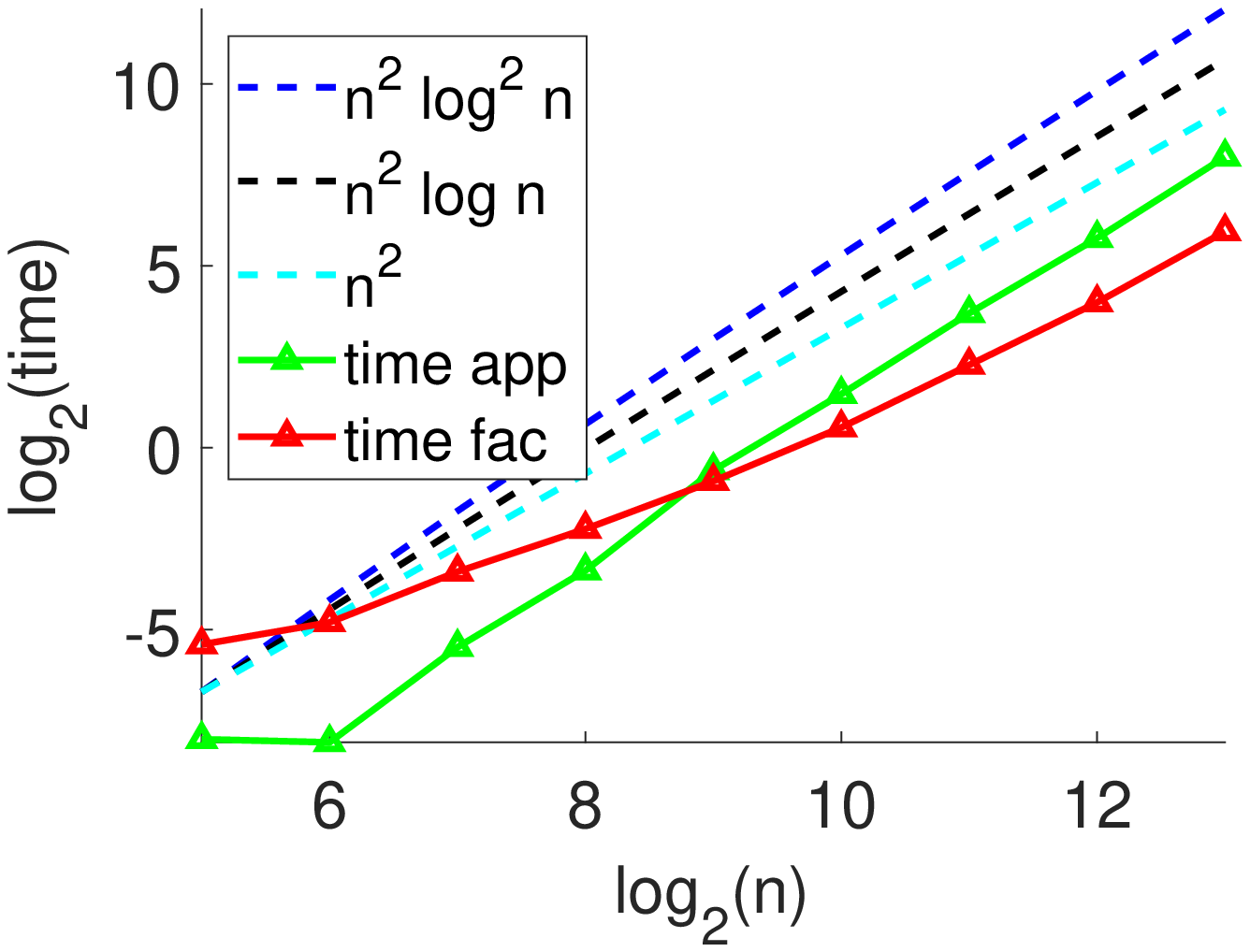}&
      \includegraphics[height=1.45in]{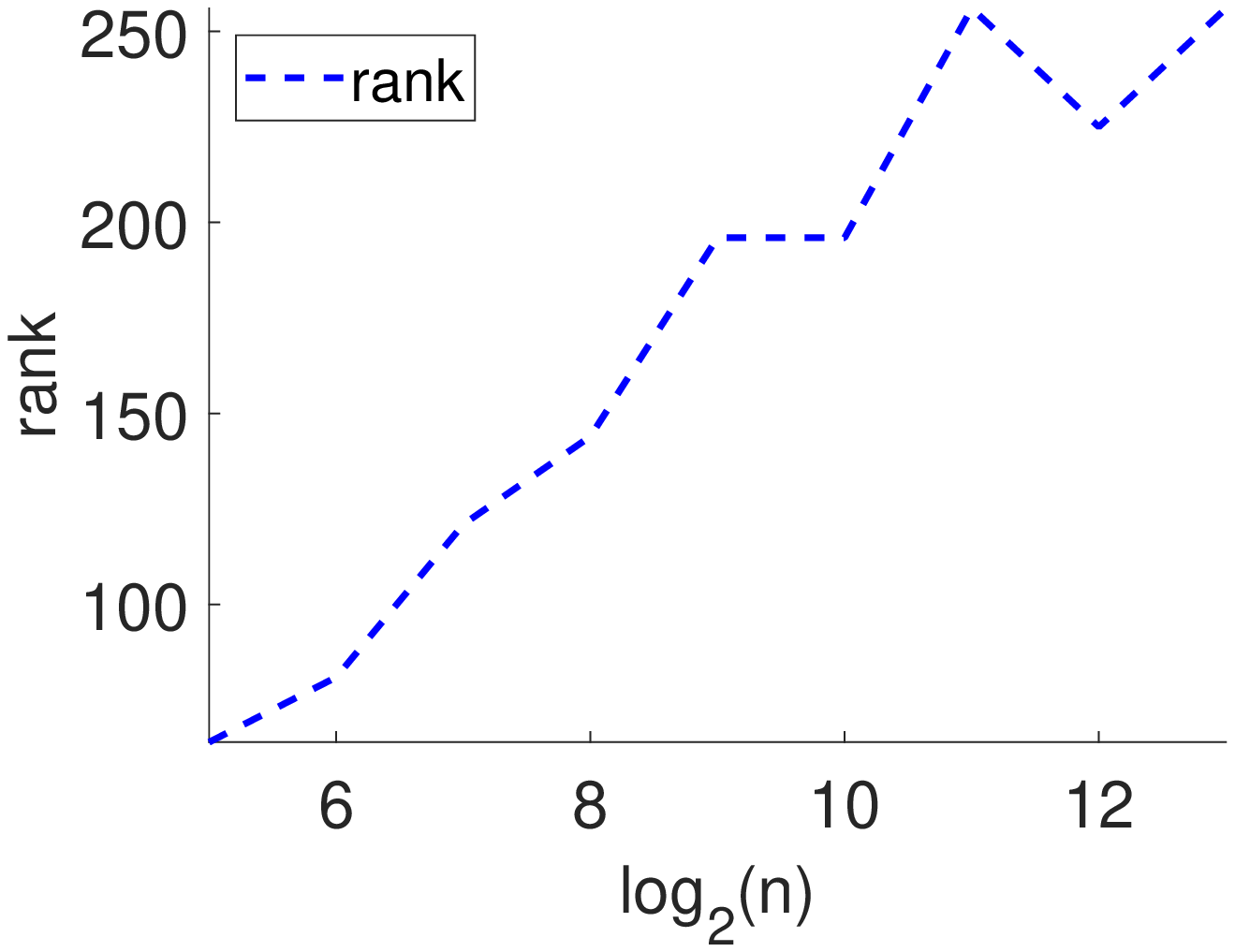}&
      \includegraphics[height=1.45in]{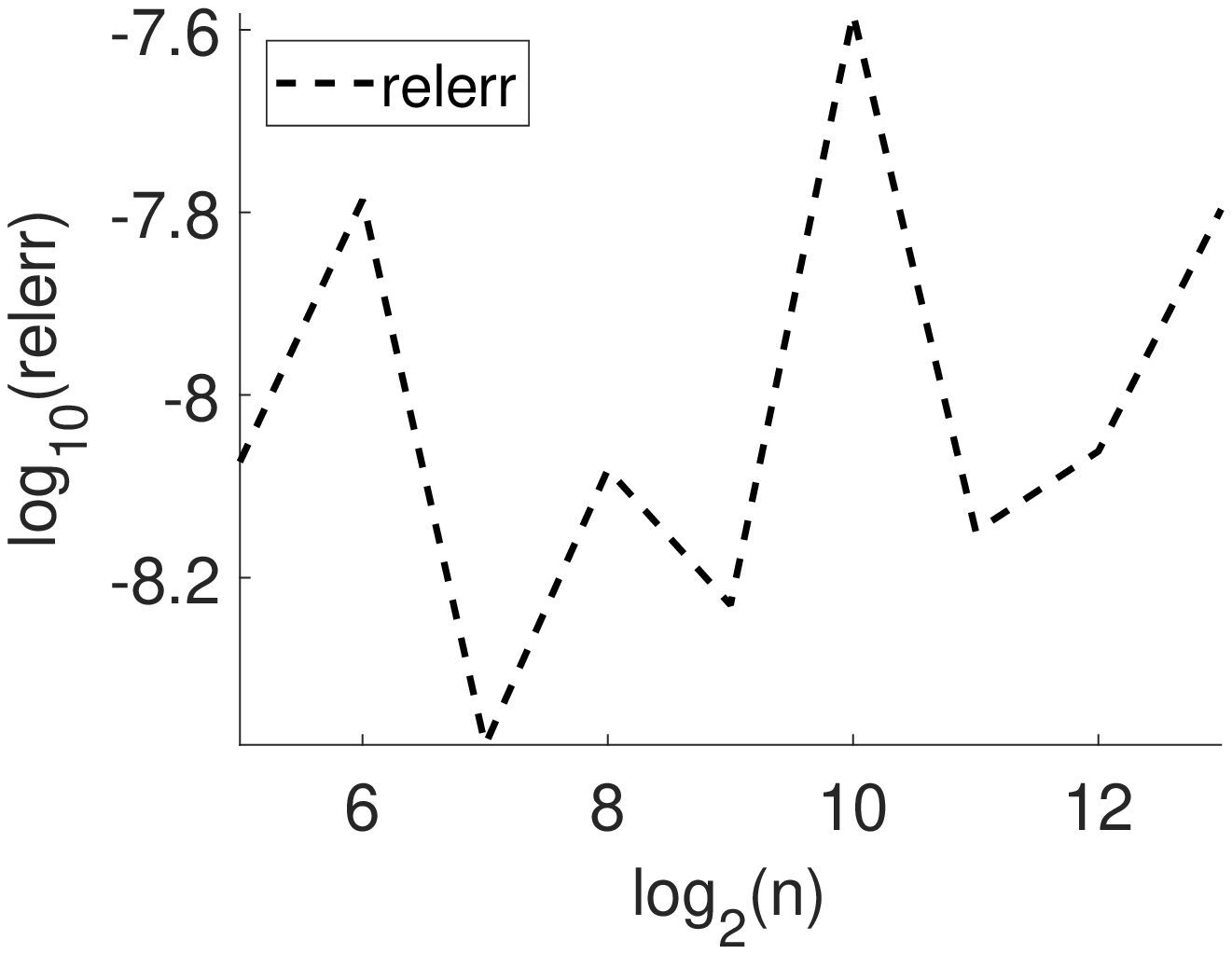}\\
      \\
      \includegraphics[height=1.45in]{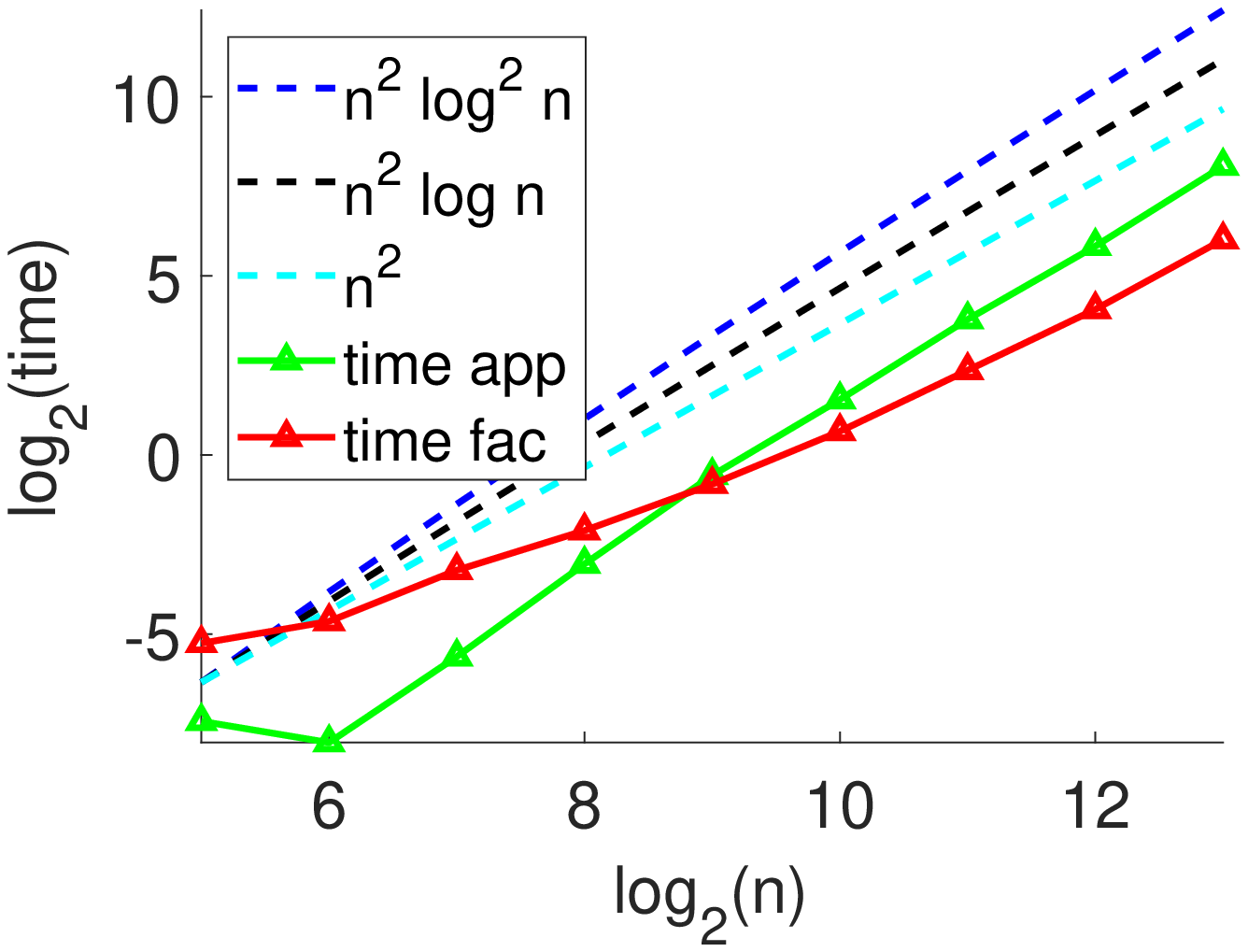}&
      \includegraphics[height=1.45in]{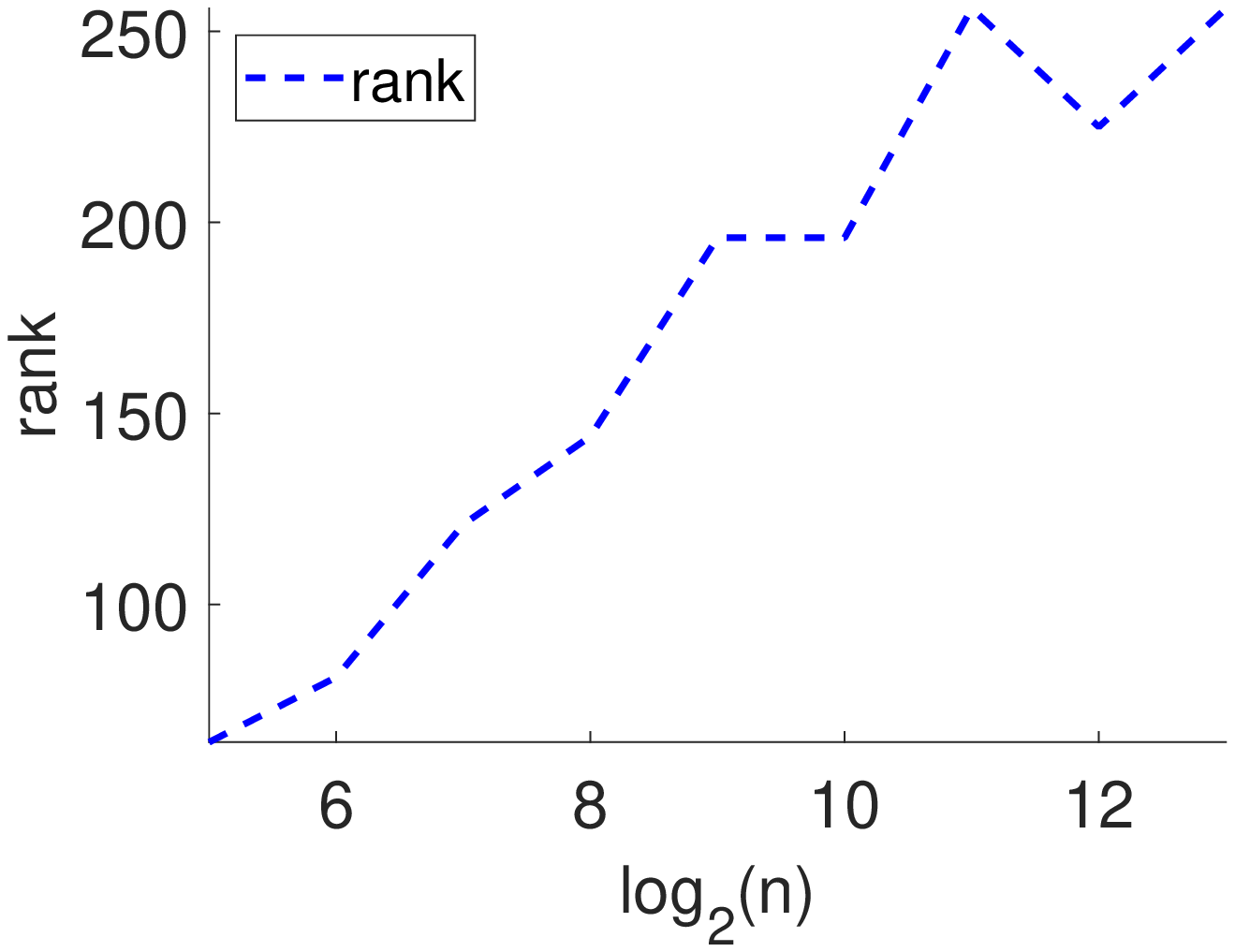}&
      \includegraphics[height=1.45in]{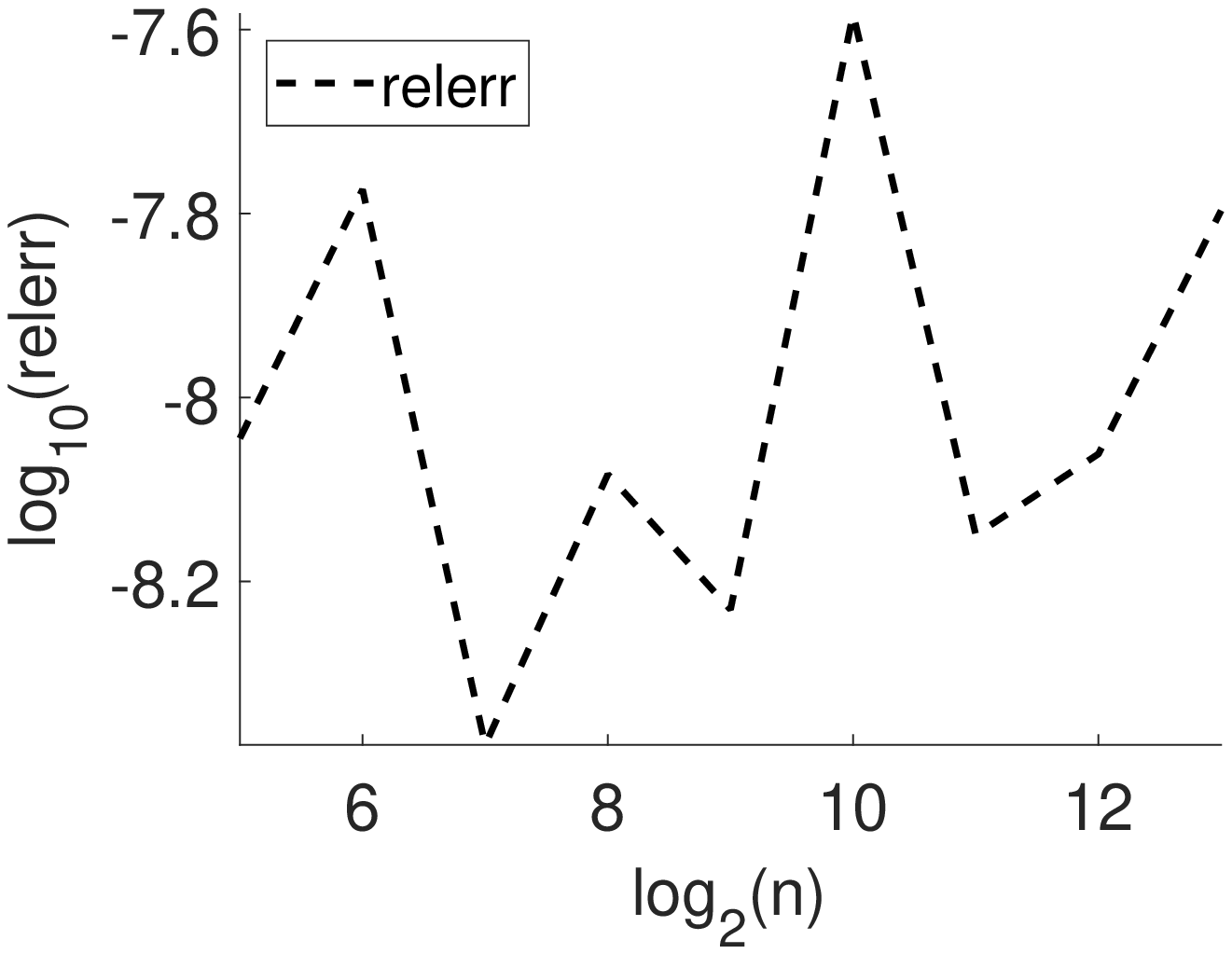}\\
      \\
      \includegraphics[height=1.45in]{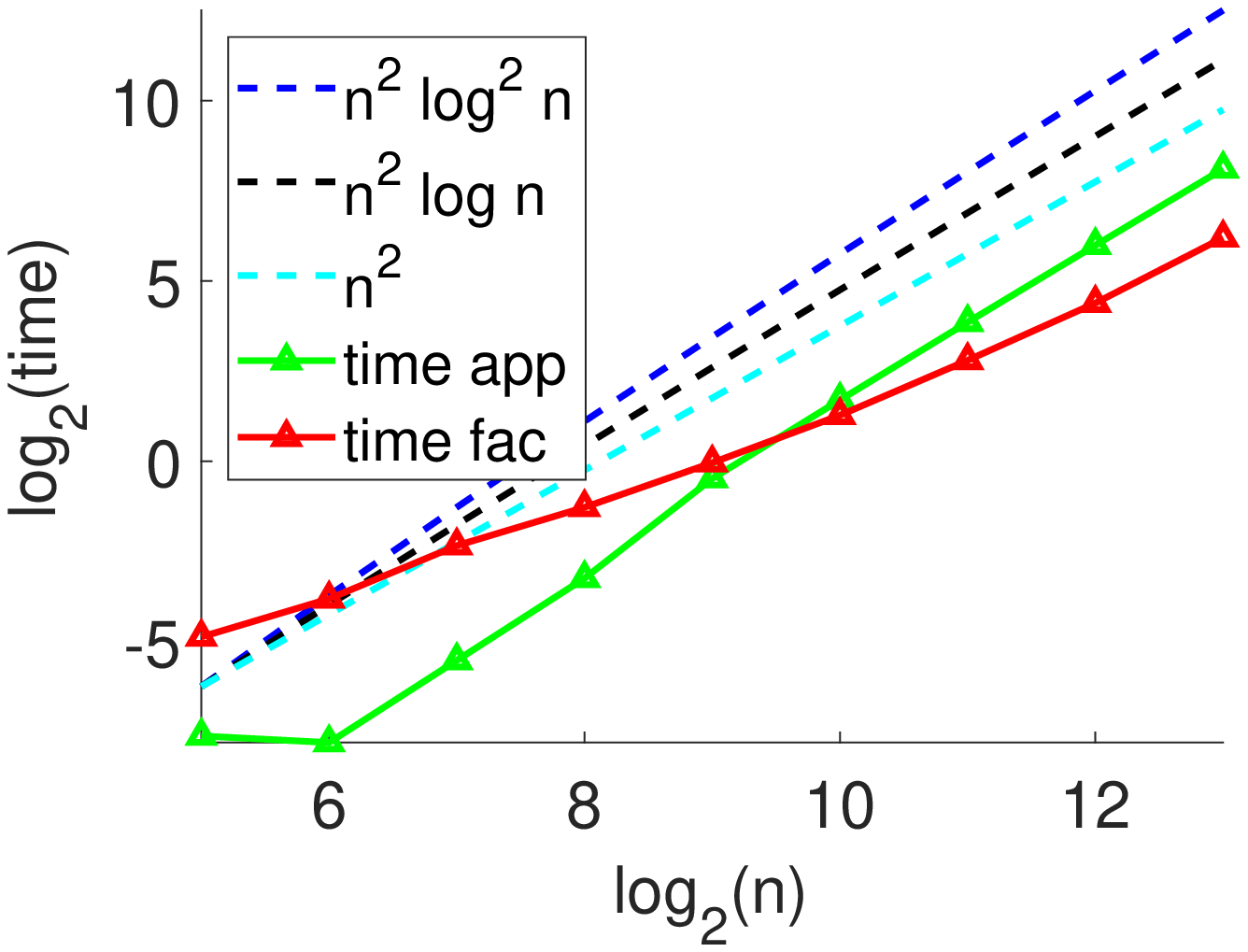}&
      \includegraphics[height=1.45in]{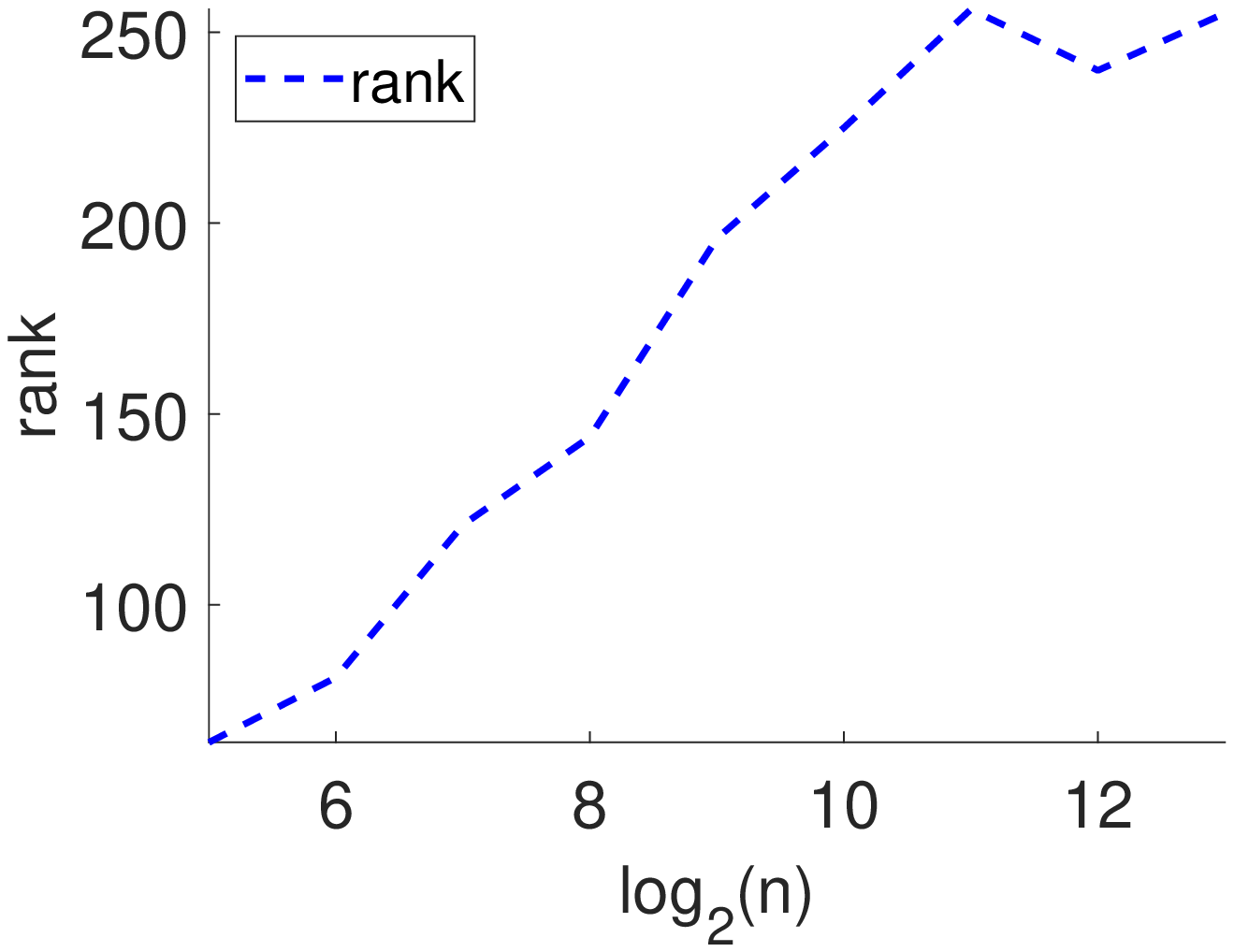}&
      \includegraphics[height=1.45in]{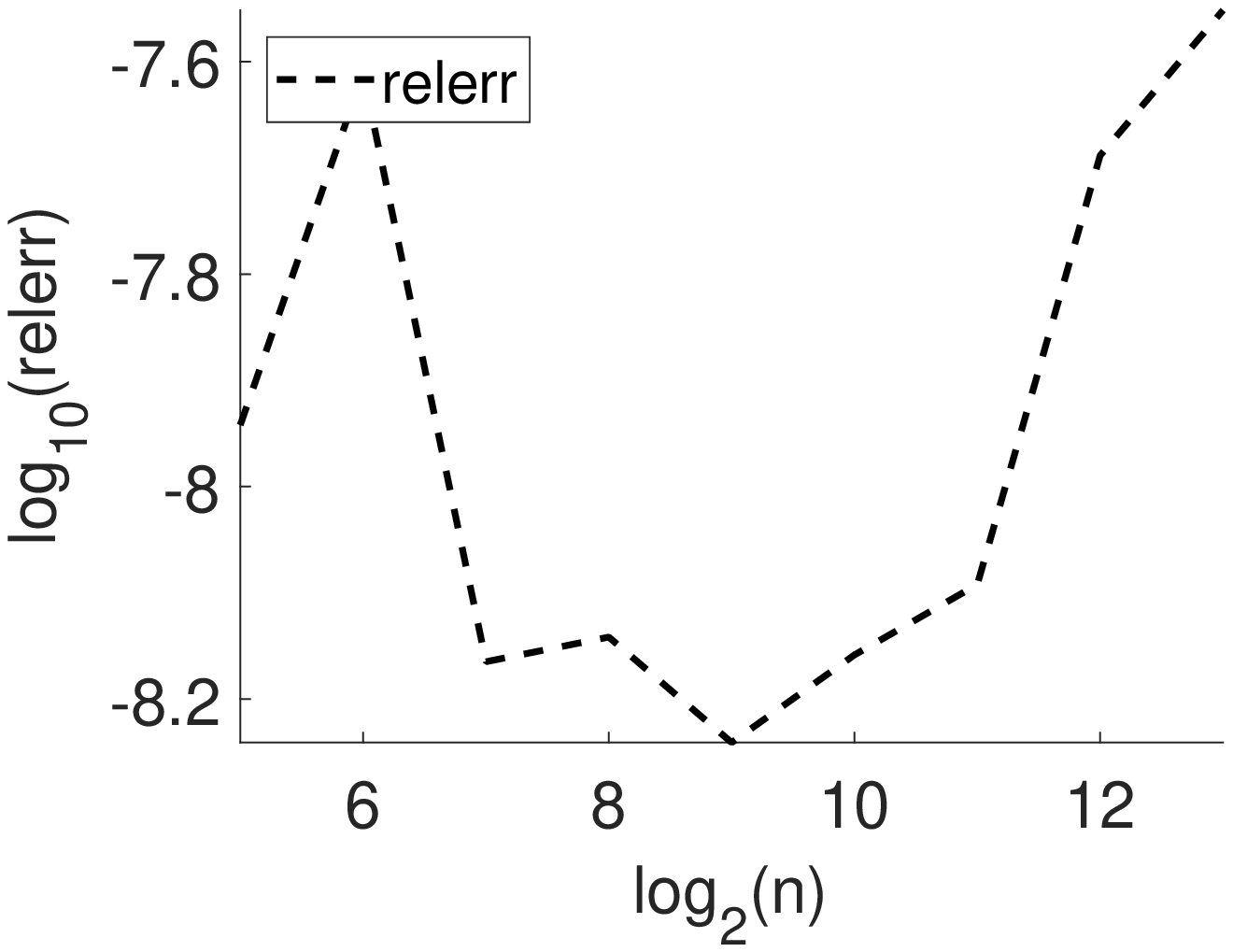}\\
    \end{tabular}
  \end{center}
\caption{Numerical results for the 2D uniform forward transform (the first row), the 2D uniform inverse transform (the second row), and the 2D nonuniform forward transform. The running time  {{(in second)}}, the numerical rank of the low-rank matrix in \eqref{eqn:Blr}, and the relative error of the fast algorithm compared to the exact summation are visualized from left to right columns.}
\label{fig:2DwoD}
\end{figure}

    \begin{figure}[ht!]
  \begin{center}
    \begin{tabular}{ccc}
      \includegraphics[height=1.45in]{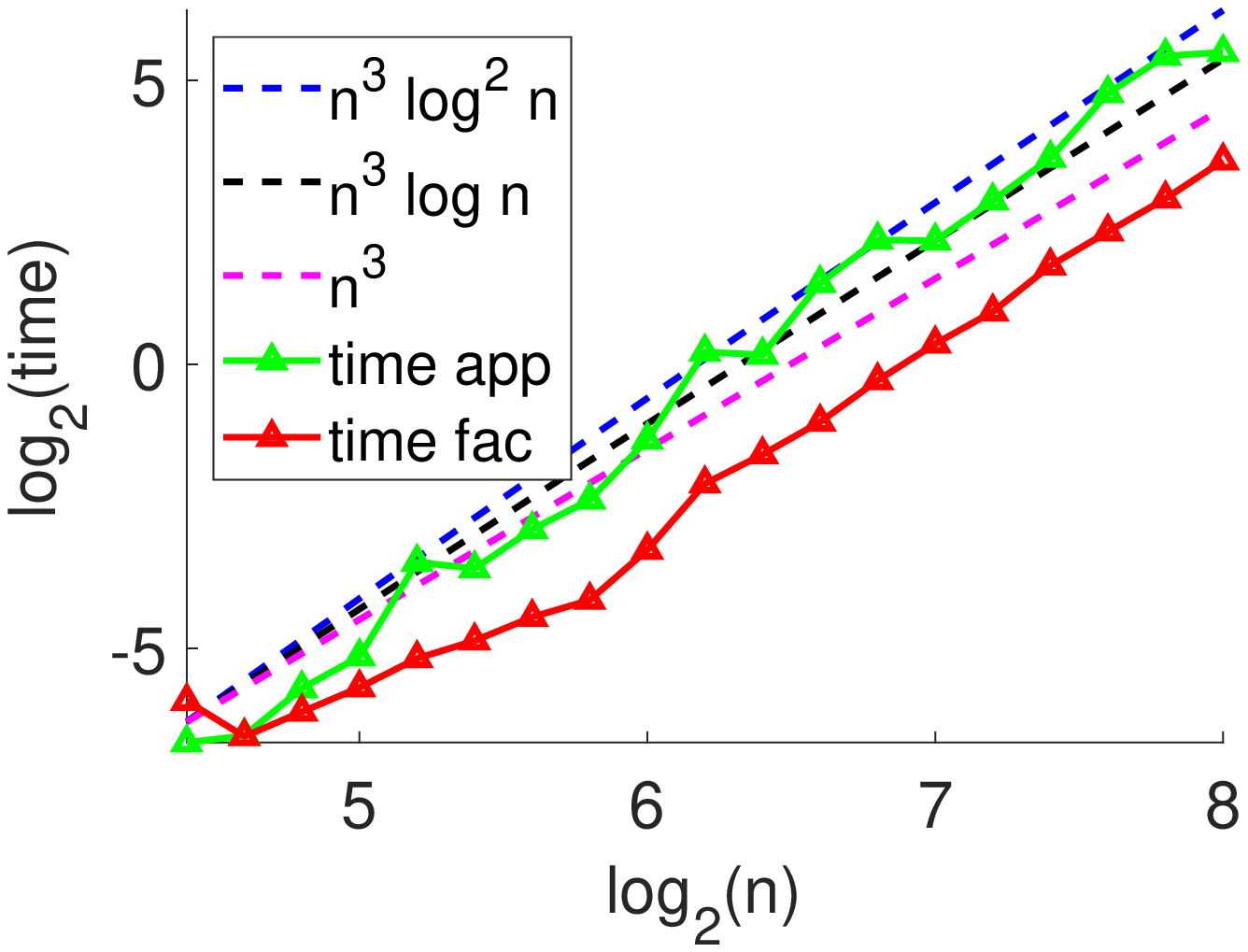}&
      \includegraphics[height=1.45in]{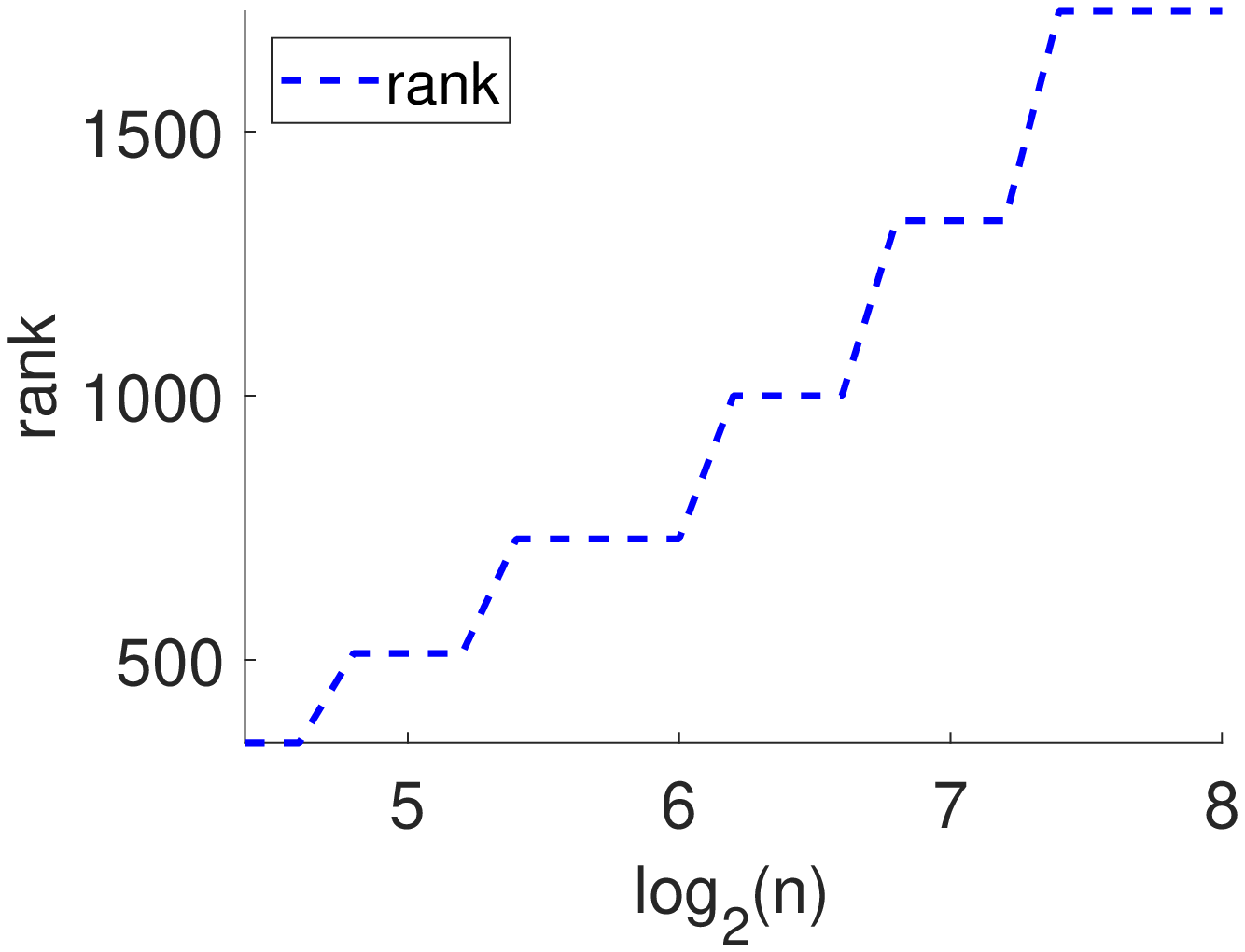}&
      \includegraphics[height=1.45in]{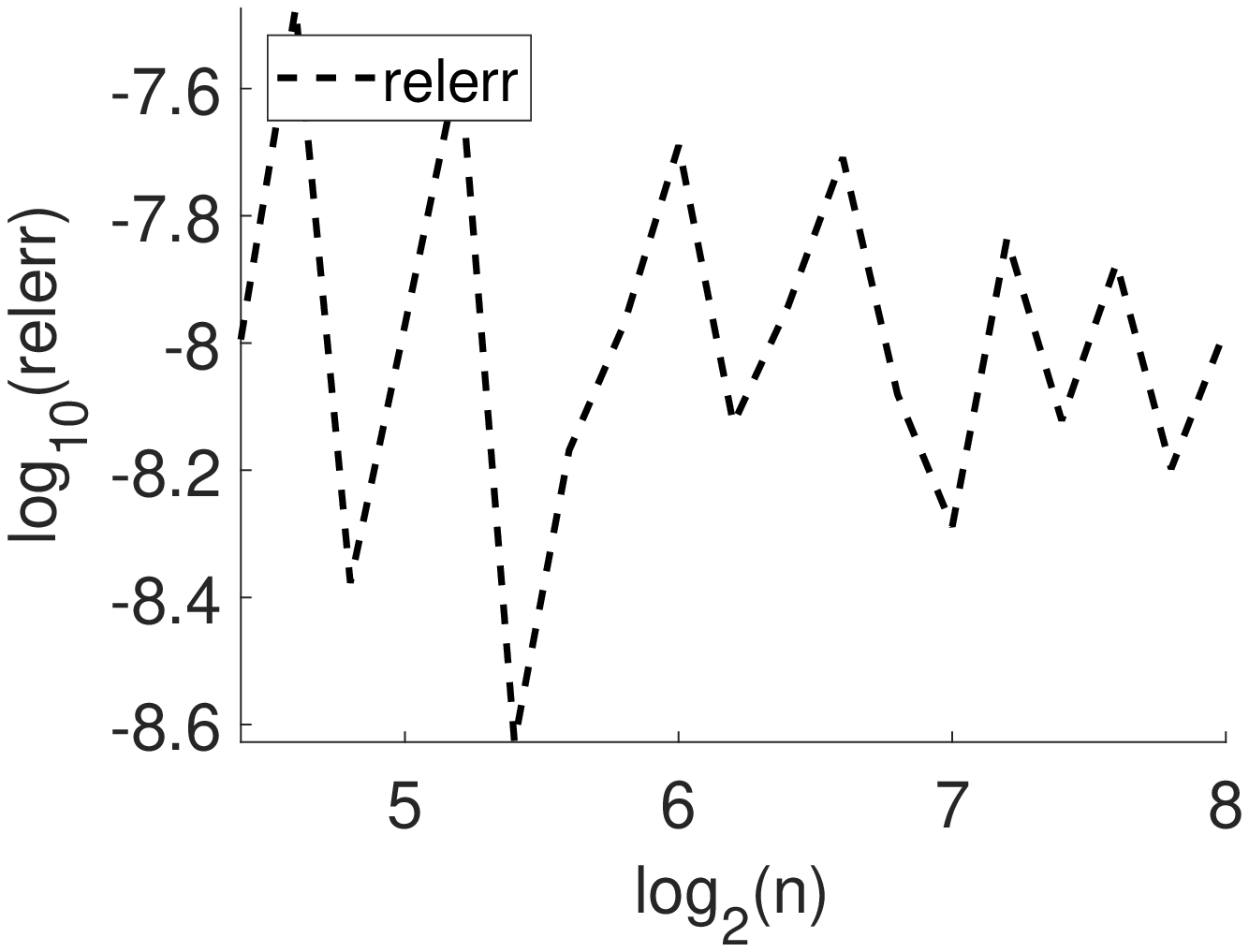}\\
      \\
      \includegraphics[height=1.45in]{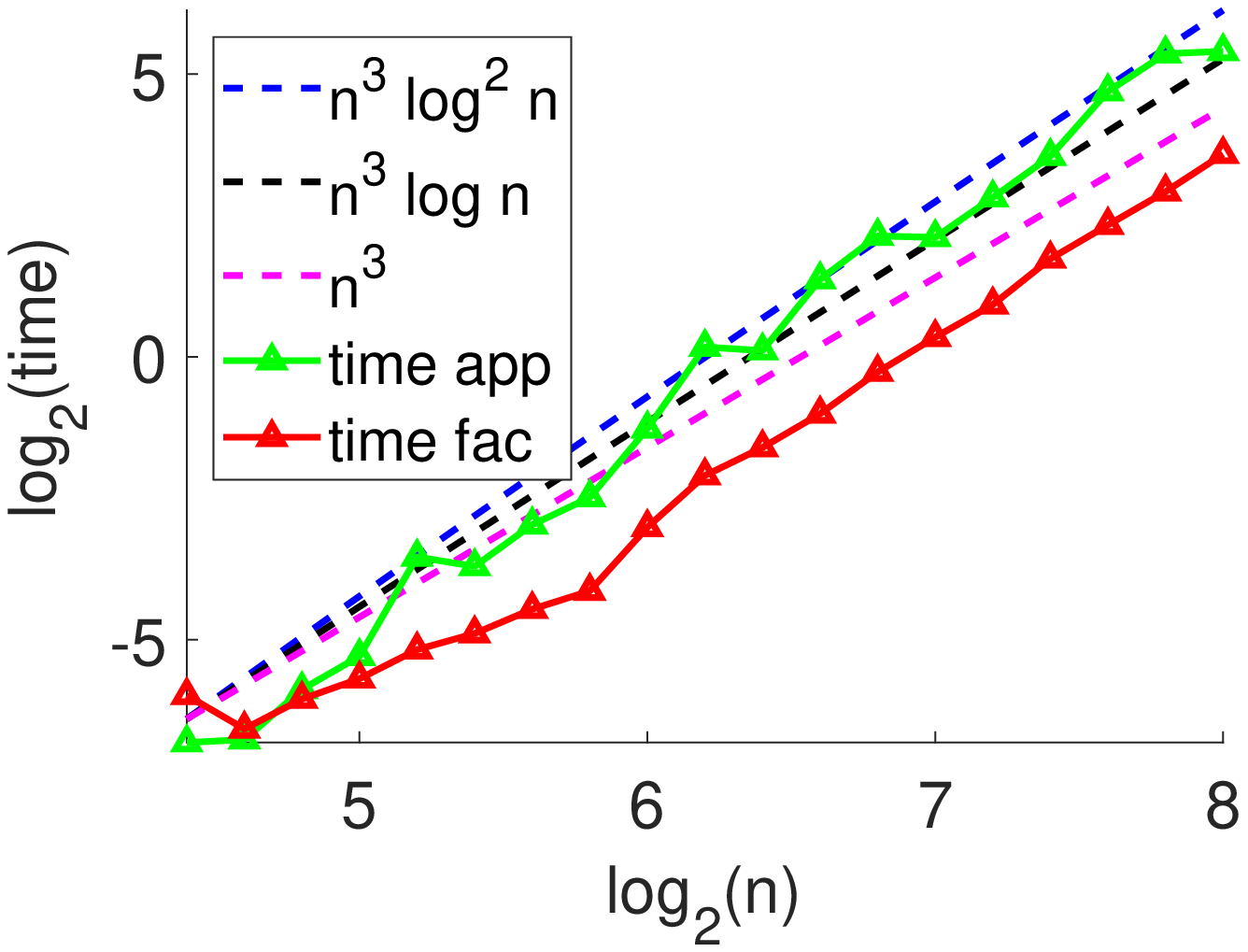}&
      \includegraphics[height=1.45in]{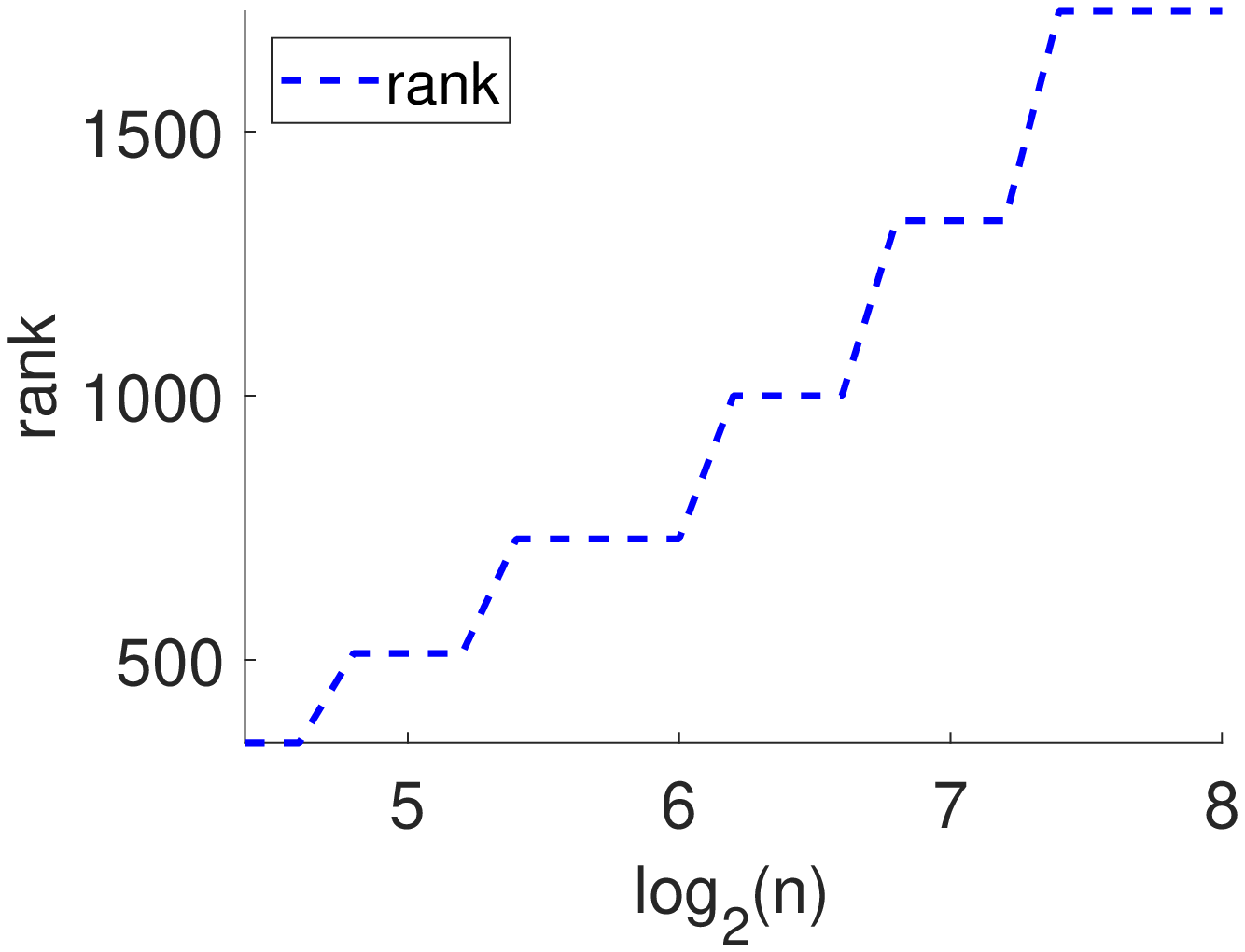}&
      \includegraphics[height=1.45in]{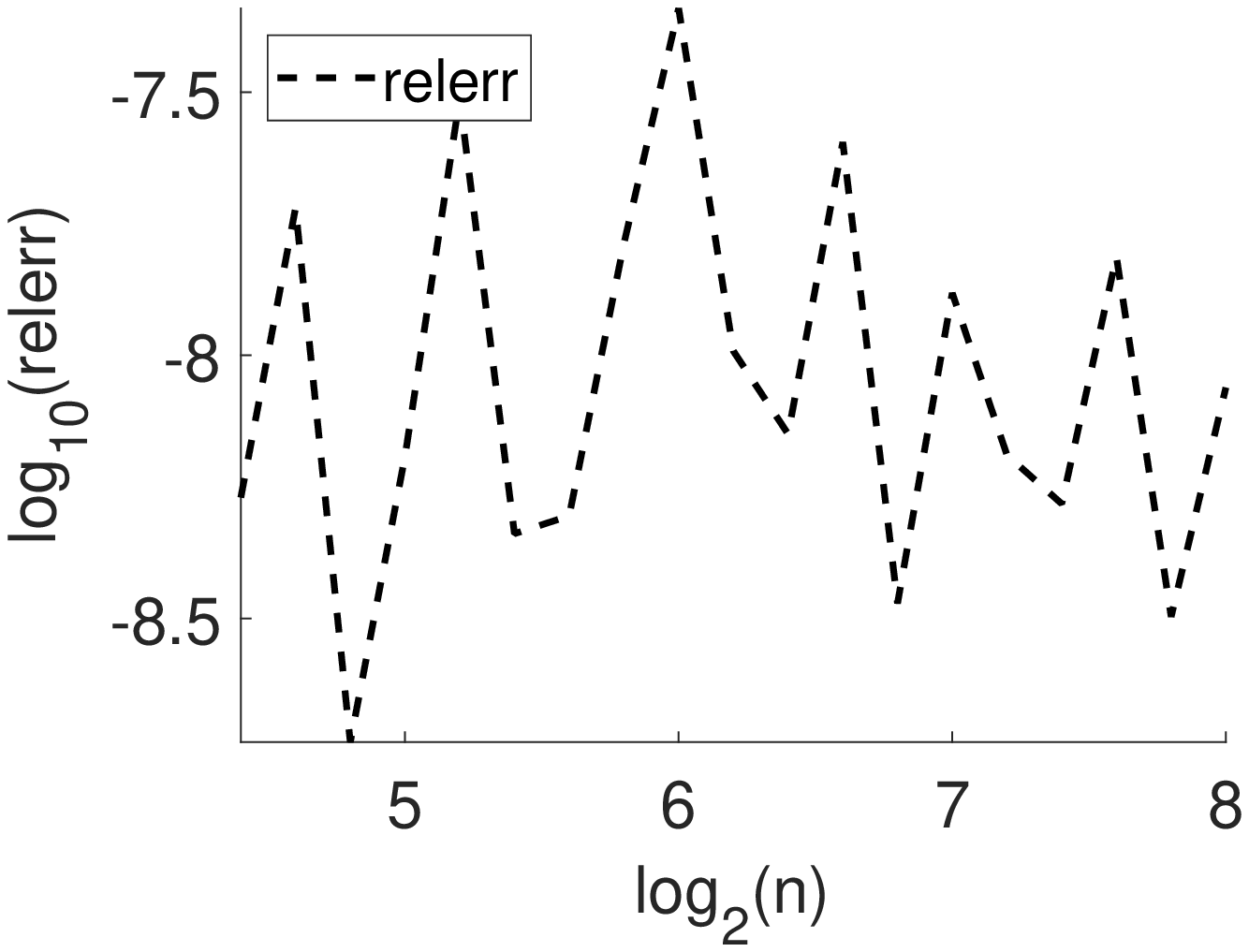}\\
      \\
      \includegraphics[height=1.45in]{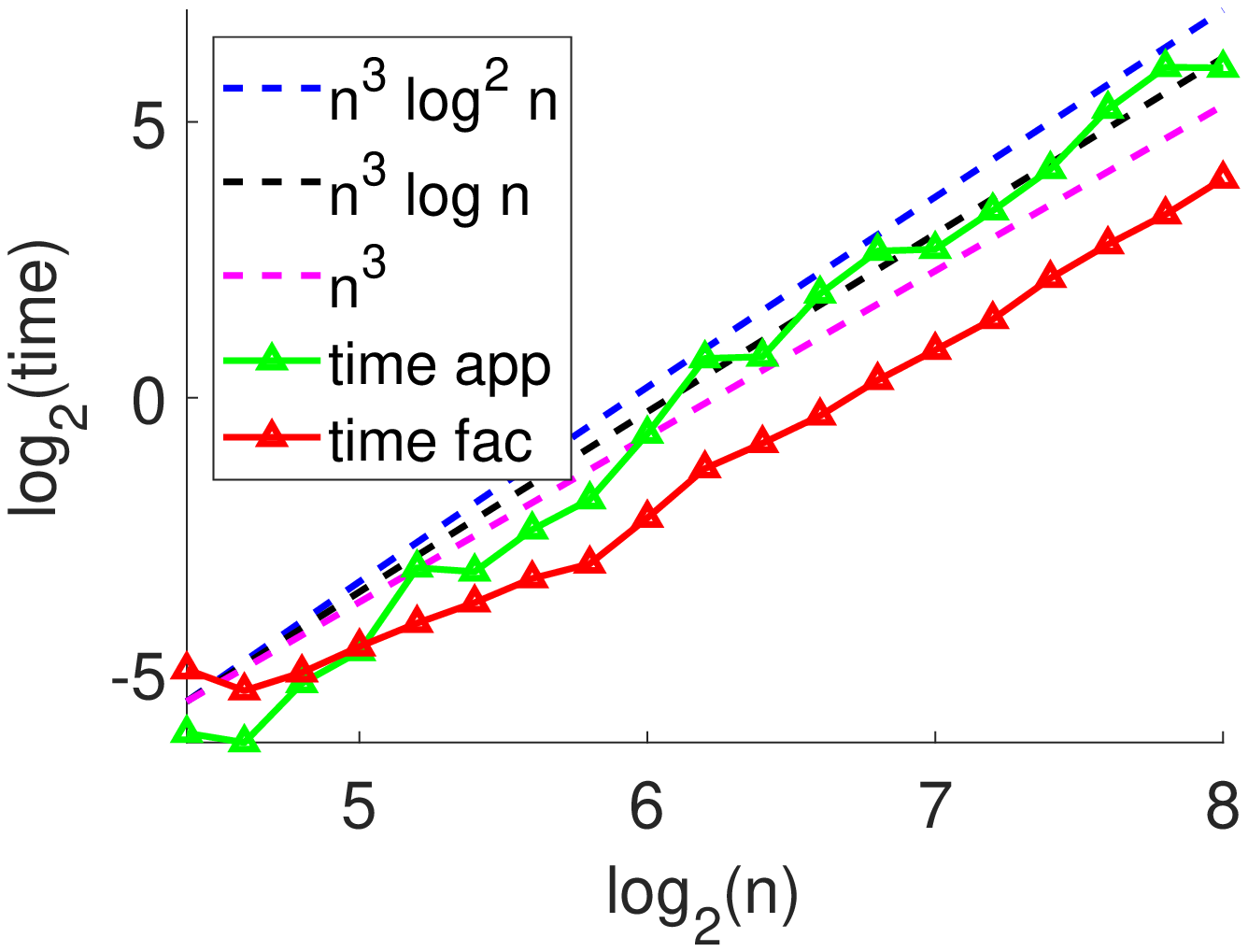}&
      \includegraphics[height=1.45in]{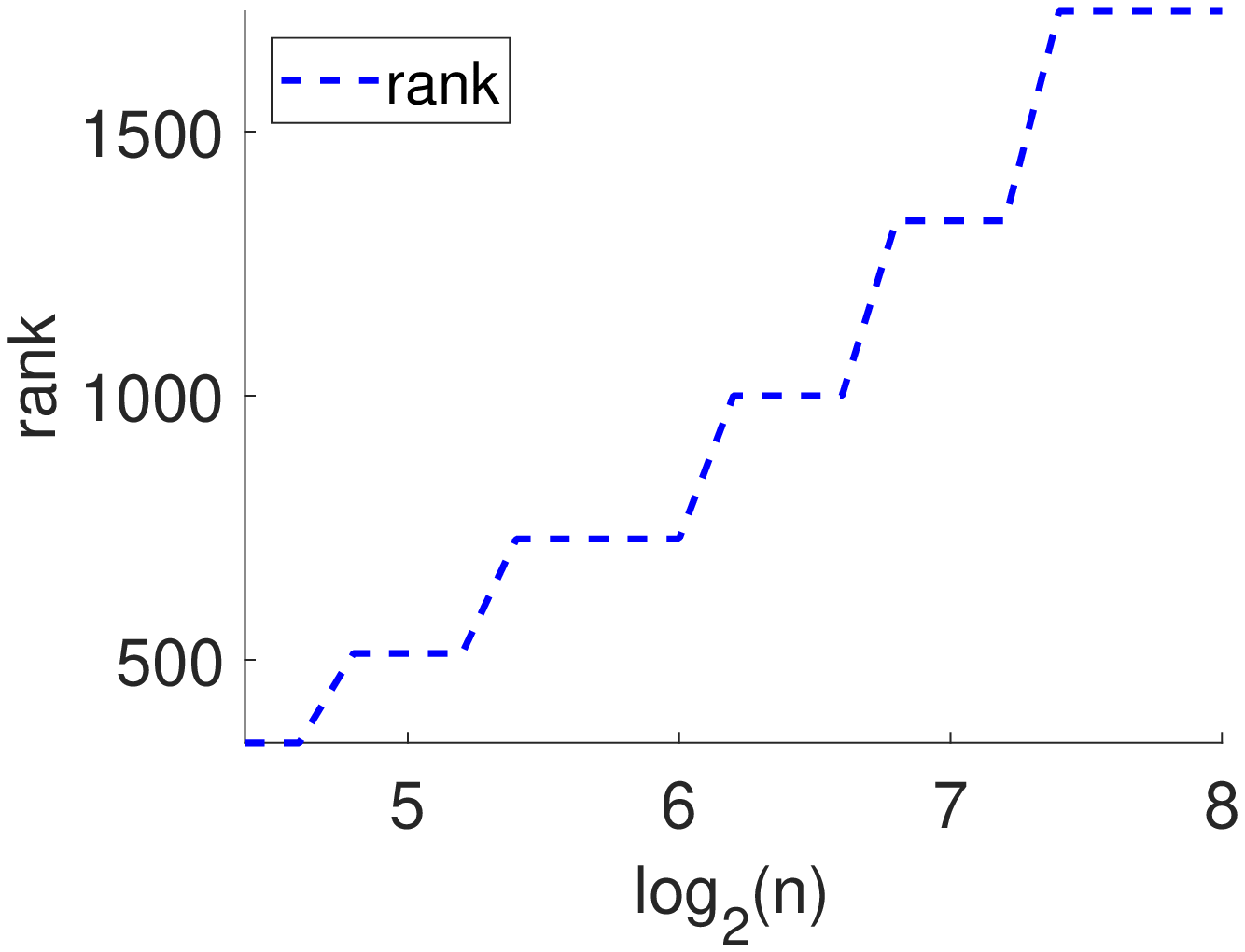}&
      \includegraphics[height=1.45in]{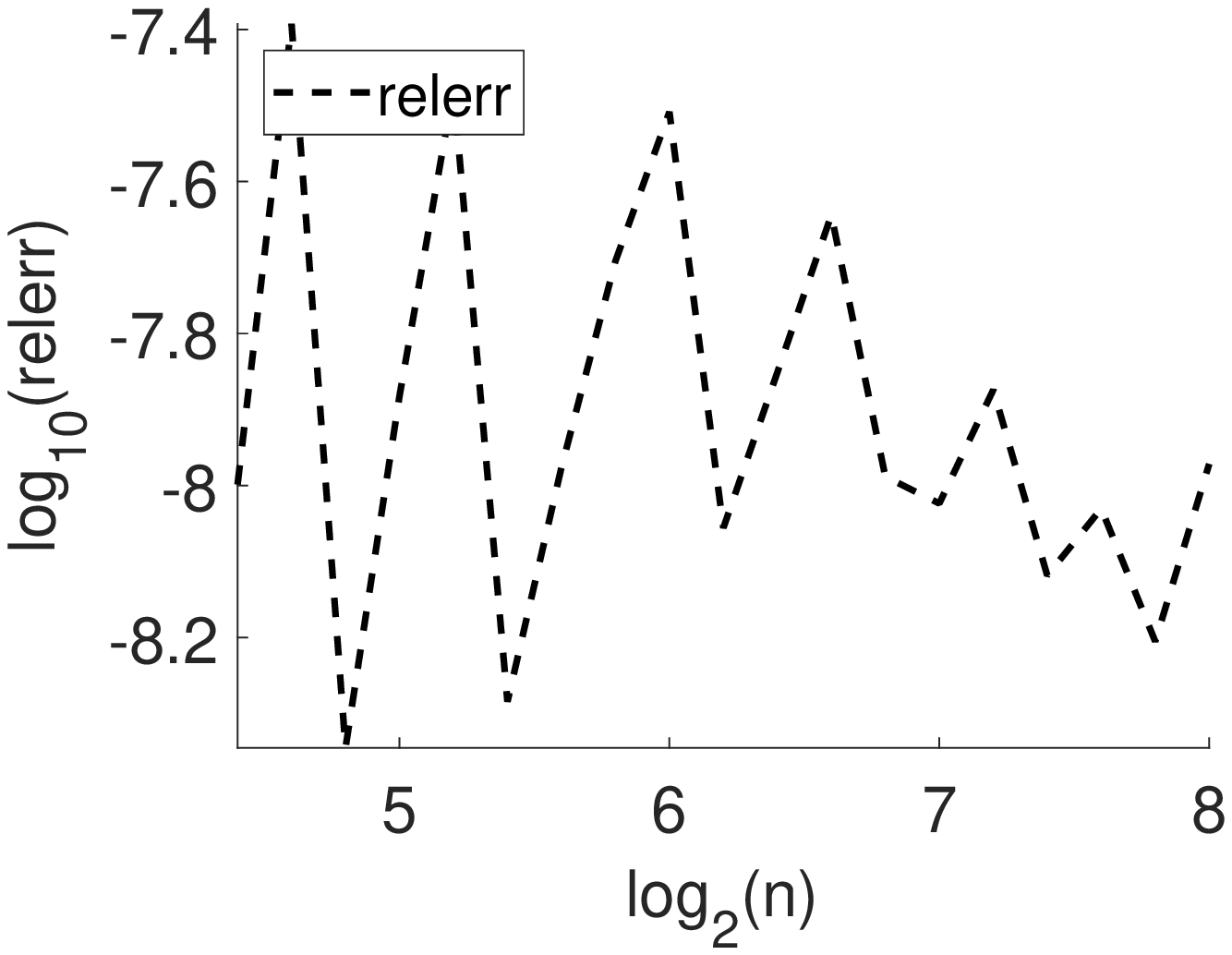}\\
    \end{tabular}
  \end{center}
\caption{Numerical results for the 3D uniform forward transform (the first row), the 3D uniform inverse transform (the second row), and the 3D nonuniform forward transform. The running time  {{(in second)}}, the numerical rank of the low-rank matrix in \eqref{eqn:Blr2}, and the relative error of the fast algorithm compared to the exact summation are visualized from left to right columns.}
\label{fig:3DwoD}
\end{figure}

In our numerical results, $n$ denotes the number of grid points per dimension; ``fac" and ``app" stand for the factorization time for setting up the algorithm and the application time for applying the algorithm, respectively; ``relerr" means the relative error of the fast matvec. The rank parameter $r$ ({{the upper bound of the rank for each dimension of a low rank factorization}}) and the accuracy parameter $\epsilon$ in all low-rank factorization algorithms are set to $\ceil[\big]{2\log_2(N)}$ and $1e-8$, respectively, where $\ceil[\big]{\cdot}$ is the ceiling function. To make fair comparisons, the oversampling parameter $q$ in Algorithm \ref{algo} for the RS method is set such that $q\times r$ is approximately equal to the number of piecewise Chebyshev grids in  (\ref{precomp:tnodes}) used in the CHEB method. We perform the same experiment $10$ times and summarize the average statistics in the following figures. {{All timings are measured in seconds.}}

{{
\subsection{Stability}
\label{sec:STA}

%
In the case of ``uniform" transform, the transformation matrix is orthonormal after weighting (Eq. \ref{eq:oneJ}) and hence the condition number is $1$. Hence, the forward and inerse transforms are well-conditioned. The proposed fast algorithm explores low-rank approximation to special structures in the Jacobi polynomial matrix $\mathcal{J}_{n,x}^{(a,b)}$ and may result in perturbation error. To ease the concern of  numerical stability, we evaluate the following quantity for random vectors $v$ to measure the stability of the propose fast forward and inverse transforms:
\begin{equation}
\label{eq:stability}
\| (W \mathcal{J}_{n,x}^{(a,b)})^{-1}W \mathcal{J}_{n,x}^{(a,b)} v - v \|/\|v\|,
\end{equation}
where $(W \mathcal{J}_{n,x}^{(a,b)})^{-1} = (W \mathcal{J}_{n,x}^{(a,b)})^{T}$ due to the orthonormal property. Table \ref{tb:stability} summarizes the above quantity for different $a$'s and $b$'s in different dimensions and shows that the proposed fast algorithm is stable.}}

\subsection{Comparison of CHEB and RS methods}
\label{sec:CDLFA}

\subsubsection{Comparisons in the one-dimensional case for $a$ and $b$ in $(-1,1)$}


{{First, we provide comparisons between CHEB and RS algorithms when $a$ and $b$ are set to be different values in the range $(-1,1)$. The numerical results are summarized in Table \ref{table1} and Table \ref{table2}. The numerical results show that both CHEB and RS achieve desired accuracy when $a$ and $b$ are not close to $-1$; RS works but CHEB fails when $a$ and $b$ are close to $-1$. Hence, the newly proposed RS algorithm works in a larger range of $a$ and $b$ than the CHEB method in \cite{Jacobi}.}}

\subsubsection{Comparisons in the three-dimensional case for $a$ and $b$ in $(-1,1)$}

Second, we provide performance comparisons of CHEB algorithm and RS algorithm in the three-dimensional case with $a=b=0.40$ to demonstrate that the RS is more efficient than the CHEB algorithm. In the one and two-dimensional cases, RS is also faster than CHEB but the speed-up is less obvious. 

As we can see from the results summarized in Figure \ref{fig:comp3}, the RS method provides a more compact matrix compression while keeping the compression accuracy competitive to that of the CHEB algorithm. The more compact compression results in faster set-up and application time for the Jacobi transform for all problem sizes.

\subsection{Performance of fast multi-dimensional transforms}
\label{sec:PFMT}
Finally, we will present examples for two and three-dimensional Jacobi polynomial transforms using the RS algorithm. In either the two or three-dimensional case, there are numerical results for one uniform forward transform, one uniform inverse transform, and one nonuniform forward transform. $a$ and $b$ are all set to $0.40$ in these examples.

{{Figure \ref{fig:2DwoD} and \ref{fig:3DwoD} summarize the numerical results in 2D and 3D, respectively.}} In cases, no matter forward or inverse transform, uniform or nonuniform transform, 2D or 3D, {{the factorization and application time of our algorithm scales like $O(n^d \log n)$ or $O(n^d\log^2 n)$, where $d$ is the dimension. The numerical rank gradually increases like $O(\log n)$ as the problem size increases.}} The relative error of the fast algorithm is in line with the desired accuracy in the low-rank factorization.

\section{Conclusion}
\label{sec:conclusion}

This paper proposed a fast algorithm for multi-dimensional Jacobi polynomial transforms based on the observation that the solution of Jacobi's differential equation can be represented via non-oscillatory phase and amplitude functions. In particular, the transformation matrix corresponding to the discrete transform is a Hadamard product of a numerically low-rank matrix and a multi-dimensional discrete Fourier transform matrix. After constructing a low-rank approximation to the numerical low-rank matrix, the application of the Hadamard product can be carried out via $O(r^d)$ fast Fourier transforms, where $r=O(\dfrac{\log n}{\log \log n})$, $n$ is the grid size per dimension, and $d$ is the dimension, resulting in a nearly optimal algorithm to compute the multi-dimensional Jacobi polynomial transform. 

We proposed to apply the randomized SVD to construct low-rank factorizations, resulting in a faster Jacobi transform in high-dimensional spaces. Moreover, numerical experiments show that the new fast transform works for a larger class of Jacobi polynomials with parameter $a$ and $b$
in the interval $(-1,1)$ than the one-dimensional algorithm in \cite{Jacobi}.

For other values of $a$ and $b$ outside $(-1,1)$, the Jacobi transformation matrix is no longer purely oscillatory and hence the proposed method is not applicable. We will tackle this issue in the future via hierarchically dividing the transformation matrix into sub-matrices that are either purely non-oscillatory, which can be handled by fast low-rank factorization, or purely oscillatory, which can be applied via the fast algorithm in this paper.

The fast inverse in the nonuniform case is still an open problem. It involves a highly ill-conditioned linear system of equations. We are working on an efficient preconditioner for this linear system and will summarize our work in a separate paper.

{\bf Acknowledgments.} H. Y. was partially supported by Grant R-146-000-251-133 in the Department of Mathematics at the National University of Singapore, by the Ministry of Education in Singapore under the grant MOE2018-T2-2-147, and the start-up grant of the Department of Mathematics at Purdue University.

\bibliographystyle{unsrt} 
\bibliography{ref}

\end{document}